\documentclass{article}
\usepackage{amsmath, amssymb}
\newtheorem{theorem}{Theorem}

\newtheorem{definition}[theorem]{Definition}
\newtheorem{example}[theorem]{Example}
\newtheorem{remark}[theorem]{Remark}
\newtheorem{corollary}[theorem]{Corollary}

\newcommand{\examp}[1]
  {\begin{example} {\rm #1} \end{example}}
\newcommand{\rem}[1]
  {\begin{remark} {\rm #1} \end{remark}}

\def\QED{\quad\blackslug\lower 8.5pt\null}

\newcommand{\crazy}[2]{\displaystyle{\mathop{#1}_{#2}}
\vphantom{\displaystyle{#1}}}

\begin{document}

\begin{center}
{\Large\bf A CLASSIFICATION AND EXAMPLES}

 \vspace*{2mm}

{\Large\bf OF FOUR-DIMENSIONAL NONISOCLINIC  }

\vspace*{2mm}

{\Large\bf   THREE-WEBS}

\vspace*{4mm}

{\large Vladislav V. Goldberg}
\end{center}

 Abstract. {\footnotesize  A classification and examples  of 4-dimensional
nonisoclinic 3-webs
of codimension 2 are given. The examples considered prove the existence
for many classes of webs for which the general existence
theorems are not proved yet.}

 \vspace*{8mm}

\setcounter{section}{-1}

\section{Introduction} In the 1980s while studying rank problems
for webs (see the papers [G 83, 92], the survey paper [AG 00]
and the monograph [G 88], Ch. 8), the author has constructed three
examples of exceptional  4-webs $W (4, 2, 2)$ of maximum 2-rank on a
4-dimensional manifold  $X^4$ (see [G 85, 86, 87], [AG 00], and
the books [G 88], [AS 92], and [AG 96]). They are exceptional
since  they are of maximum rank but not  algebraizable.
 H\'{e}naut [H 98]  named these webs after the author
and denoted them by $\mathcal{G}_1 (4, 2, 2),
\mathcal{G}_2 (4, 2, 2)$, and $\mathcal{G}_3 (4, 2, 2)$.

When the author was constructing these examples of 4-webs, he
 considered numerous examples of 3-webs $W (3, 2, 2)$ on $X^4$ and
proved that  3 of them can be  expanded  to exceptional
 4-webs $\mathcal{G}_1 (4, 2, 2), \mathcal{G}_2 (4, 2, 2)$,
 and $\mathcal{G}_3 (4, 2, 2)$
 (see [G 85, 86, 87, 88]).

  However, it turns out that many  examples
of 3-webs $W (3, 2, 2)$ that the author has constructed in that research
 are useful when one studies different classes of multidimensional
3-webs (isoclinic, hexagonal, transversally
geodesic, algebraizable, Bol's webs, etc.). Some of these examples
of  webs $W (3, 2, 2)$ were published in
the author paper [G 92] and the book [G 88], and some were
included in problem sections of the book [AS 92].

The author decided to present some of these examples both
published and unpublished following some classification for them.
It is worth to consider these examples in order to provide an up-to-date
characterization of webs by indicating to which classes
 a web belongs.

For isoclinic webs this was done in [G 99]. In the current
paper we present a classification and
examples of nonisoclinic webs $W (3, 2, 2)$.

The examples mentioned above prove  the existence
 of many classes of webs for which the general existence theorems
are not proved yet.

 \section{The transversal distribution of a web {\protect\boldmath\(W (3, 2, 2)\)
\protect\unboldmath}}

\textbf{1.} The leaves of the foliation $\lambda_\xi, \; \xi =
1, 2, 3$,  of a web $W (3, 2, 2)$ are determined by the equations
$\crazy{\omega}{\xi}^i = 0, \; i = 1, 2$,   where
\begin{equation}\label{eq:1}
\crazy{\omega}{1}^i + \crazy{\omega}{2}^i + \crazy{\omega}{3}^i = 0
\end{equation}
(see, for example, [G 88], Section {\bf 8.1} or
[AS 92], Section  {\bf 1.3}).

The structure equations of such a web can be written in the form
\begin{equation}\label{eq:2}
\renewcommand{\arraystretch}{1.3}
 \left\{
\begin{array}{ll}
   d \crazy{\omega}{1}^i =
\crazy{\omega}{1}^j \wedge \omega_j^i +  a_j \crazy{\omega}{1}^j
\wedge \crazy{\omega}{1}^i, \\  d\crazy{\omega}{2}^i =
\crazy{\omega}{2}^j \wedge \omega_j^j - a_j \crazy{\omega}{2}^j
\wedge \crazy{\omega}{2}^i.
\end{array}
      \right.
\renewcommand{\arraystretch}{1}
 \end{equation}
 The differential prolongations of equations (2) are
\begin{equation}\label{eq:3}
  d\omega_j^i - \omega_j^k \wedge
\omega_k^i =   b_{jkl}^i \crazy{\omega}{1}^k \wedge
\crazy{\omega}{2}^l,
 \end{equation}
\begin{equation}\label{eq:4}
da_i - a_j \omega_i^j = p_{ij} \crazy{\omega}{1}^j  +  q_{ij}
\crazy{\omega}{2}^j,
\end{equation}
 where
 \begin{equation}\label{eq:5}
b^i_{[j|l|k]} = \delta^i_{[k} p_{j]l}, \;\; b^i_{[jk]l} =
\delta^i_{[k} q_{j]l}
\end{equation}
(see [G 88], Sections {\bf 8.1} and {\bf 8.4} or [AS 92], Section
{\bf 3.2}).
The quantities
\begin{equation}\label{eq:6}
 a_{jk}^i = a_{[j} \delta_{k]}^i
 \end{equation}
 and $b^i_{jkl}$ are the {\em
torsion and curvature tensors} of a three-web  $W (3,
2, 2)$. Note that for webs $W (3, 2, 2)$ the torsion tensor
$a^i_{jk}$ always has structure (6), where  $a = \{a_1, a_2\}$ is
a covector. If $a = 0$, then a web $W (3, 2, 2)$ is {\em isoclinicly
geodesic}. In what follows in this paper, {\em we will assume that $a
\neq 0$, i.e., a web $W (3, 2, 2)$ is nonisoclinicly geodesic.}

\textbf{2}. For a web $W (3, 2, 2)$, a transversally geodesic
distribution is defined (cf. [AS 92], Section {\bf 3.1}) by the
equations
$$ \xi^2 \crazy{\omega}{1}^1 - \xi^1
\crazy{\omega}{1}^2= 0, \;\; \xi^2 \crazy{\omega}{2}^1 - \xi^1
\crazy{\omega}{2}^2 = 0.
$$
If we take $\displaystyle
\frac{\xi^1}{\xi^2} = - \frac{a_2}{a_1}$, we obtain an
invariant transversal distribution $\Delta$ defined by the
equations
\begin{equation}\label{eq:7}
 a_1 \crazy{\omega}{1}^1 + a_2 \crazy{\omega}{1}^2 = 0, \;\; a_1
 \crazy{\omega}{2}^1 + a_2 \crazy{\omega}{2}^2 = 0.
  \end{equation}

We will call the distribution $\Delta$ defined by equat\-ions (7)
the {\em transversal $a$-distribution} of a web $W (3, 2, 2)$
since it is defined by the covector $a$.
Note that for isoclinicly geodesic webs $W (3, 2, 2)$, for which
$a_1 = a_2 =0$, the transversal distribution is not defined.

The following theorem gives the conditions of integrability of the
distribution $\Delta$ (see the proofs of this and other results
of this section in [AG 98]).

\begin{theorem}
The transversal $a$-distribution $\Delta$ defined by
 equations $(7)$  is integrable if and only if the components
$a_1$ and $a_2$ of the covector $a$ and their Pfaffian derivatives
$p_{ij}$ and $q_{ij}$ satisfy the conditions
\begin{equation}\label{eq:8}
\renewcommand{\arraystretch}{1.3}
\left\{
\begin{array}{ll}
a_2^2 p_{11} -  2a_1 a_2 p_{(12)} + a_1^2 p_{22} = 0, \\
a_2^2 q_{11} -  2a_1 a_2 q_{(12)} + a_1^2 q_{22} = 0.
        \end{array}
\right.
\renewcommand{\arraystretch}{1}
\end{equation}
\end{theorem}

\textbf{3.} For a web $W (3, 2, 2)$, it is always possible
to take a specialized frame in which  there is a relation
between the components $a_1$ and $a_2$  of the covector $a$.
For example, if the transversal distribution $\Delta$ coincides
with the distribution $\crazy{\omega}{\alpha}^1 = 0$,
%or $\crazy{\omega}{\alpha}^2 = 0$
%or $\crazy{\omega}{\alpha}^1 + \crazy{\omega}{\alpha}^2 = 0$,
 then we have $a_2 = 0$.
 % or $a_1 = 0$ or $a_1 = a_2$, respectively.
In this case,
the form $\omega_2^1$
%\omega_1^2$, and  $\omega_1^1 + \omega_1^2
%- \omega_2^1 - \omega_2^2$, respectively,
is expressed in terms of
the basis forms $\crazy{\omega}{\alpha}^i, \; \alpha = 1, 2$,
i.e., in this case we have
$
\pi_2^1 = 0,
%\;\; \pi_1^2 = 0, \;\; \pi_1^1 + \pi_1^2 - \pi_2^1 - \pi_2^2 = 0,
$
where $\pi_i^j = \omega_i^j\Bigl|_{\crazy{\omega}{\alpha}^i = 0}$.

 In examples, that we are going to present in this paper, such situations
 will occur. So, we present here the conditions of integrability for the
 4 cases that include the case indicated above.

\begin{corollary}
If for a web $W (3, 2, 2)$ one of the following conditions
\begin{equation}\label{eq:9}
a_2 = 0, \;\; \pi_2^1 = 0,
\end{equation}
\begin{equation}\label{eq:10}
a_1 = 0, \;\; \pi_1^2 = 0,
\end{equation}
\begin{equation}\label{eq:11}
a_1 = a_2, \;\; \pi_1^1 + \pi_1^2 - \pi_2^1 - \pi_2^2 = 0,
\end{equation}
\begin{equation}\label{eq:12}
a_1 = - a_2, \;\; \pi_1^1 - \pi_1^2 + \pi_2^1 - \pi_2^2 = 0
\end{equation}
holds, then the $a$-distribution $\Delta$ coincides with
the distribution $\crazy{\omega}{\alpha}^1 = 0, \;
\crazy{\omega}{\alpha}^2 = \nolinebreak 0, \linebreak
   \crazy{\omega}{\alpha}^1 + \crazy{\omega}{\alpha}^2 = 0$,
 or   $\crazy{\omega}{\alpha}^1 - \crazy{\omega}{\alpha}^2 = 0$,
   respectively.
This $a$-distribution is integrable if and only if the quantities
$p_{ij}$ and $q_{ij}$ satisfy  respectively the following
conditions:
\begin{equation}\label{eq:13}
p_{22} = q_{22} = 0,
\end{equation}
\begin{equation}\label{eq:14}
p_{11} = q_{11} = 0,
\end{equation}
\begin{equation}\label{eq:15}
p_{11} - 2 p_{(12)} + p_{22} = 0, \;\;
q_{11} - 2 q_{(12)} + q_{22} = 0,
\end{equation}
\begin{equation}\label{eq:16}
p_{11} + 2 p_{(12)} + p_{22} = 0, \;\;
q_{11} + 2 q_{(12)} + q_{22} = 0.
\end{equation}
\end{corollary}

Each of relations (8), (13), (14), (15), and (16) gives two
conditions which Pfaffian derivatives $p_{ij}$ and $q_{ij}$ of the co-vector $a$ must
satisfy in order for the $a$-distribution $\Delta$
of a web $W (3, 2, 2)$ to be integrable.

\textbf{4.} Now we state the theorem giving the conditions for
the integral surfaces of the $a$-distribution $\Delta$ to be
geodesicly parallel in some affine connections.

\begin{theorem}
If on a web $W (3, 2, 2)$ the  conditions
\begin{equation}\label{eq:17}
\renewcommand{\arraystretch}{1.3}
\left\{
\begin{array}{ll}
a_{2} p_{12} - a_1 p_{22} = 0, & a_1 p_{21} -  a_{2} p_{11} = 0, \\
a_{2} q_{12} - a_1 q_{22} = 0, &  a_1 q_{21} - a_{2} q_{11} = 0
\end{array}
\right.
\renewcommand{\arraystretch}{1}
\end{equation}
hold, then the integral surfaces $V^2$ of the $a$-distribution
$\Delta$ are geodesicly parallel in any affine connection of the
bundle of affine connections defined by the forms
\begin{equation}\label{eq:18}
  \theta_v^u =  \left(
  \begin{array}{ll}
 \theta_j^i & 0 \\
               0 & \theta_j^i
               \end{array}
     \right),
               \;\;\;\; i, j = 1, 2; \;
               u, v = 1, 2, 3, 4,
\end{equation}
where
\begin{equation}\label{eq:19}
  \theta_j^i = \omega_j^i + a^i_{jk} (p \crazy{\omega}{1}^k
  + q\crazy{\omega}{2}^k)
\end{equation}
$($see {{\rm [AS 92], p. 35}}$)$.
\end{theorem}
It follows from (4) that conditions (17) are equivalent to the following
Pfaffian equation:
$$
d t + t^2 \omega_1^2 + t (\omega_1^1 -
\omega_2^2) - \omega_2^1 = 0, \; \mbox{{\rm where}}\;\;
t = \displaystyle \frac{a_2}{a_1}.
$$
Note that it follows from Theorems 1 and 3 that
if the surfaces $V^2$ of the $a$-distribution
$\Delta$ are geodesicly parallel in any affine connection of the
bundle (18)--(19), then the  $a$-distribution
$\Delta$ of  a web $W (3, 2, 2)$ is integrable.

\begin{corollary}
If   the $a$-distribution $\Delta$ of  a web $W (3, 2, 2)$
 coincides with
the distribution $\crazy{\omega}{\alpha}^1 = 0$, or
$\crazy{\omega}{\alpha}^2 = 0$, or
 $\crazy{\omega}{\alpha}^1 + \crazy{\omega}{\alpha}^2 = 0$,
 or   $\crazy{\omega}{\alpha}^1 - \crazy{\omega}{\alpha}^2 = 0$,
 then   the integral surfaces $V^2$ of the $a$-distribution
$\Delta$ are geodesicly parallel in any affine connection of the
bundle  $(18)$--$(19)$  if and only if the quantities
$p_{ij}$ and $q_{ij}$ satisfy  respectively the following
conditions:
\begin{equation}\label{eq:20}
 p_{2i} =  q_{2i} = 0,
\end{equation}
\begin{equation}\label{eq:21}
 p_{1i} =  q_{1i} = 0,
\end{equation}
\begin{equation}\label{eq:22}
 p_{1i} = p_{2i}, \;\;  q_{1i} = q_{2i},
\end{equation}
\begin{equation}\label{eq:23}
 p_{1i} = - p_{2i}, \;\;   q_{1i} = -q_{2i}.
\end{equation}
\end{corollary}

Note that conditions (20)--(23) are equivalent to
the following relations between the forms $\omega_i^j$
(cf. the remark after Theorem 3):
for these 4 possible specializations,
  $$
\renewcommand{\arraystretch}{1.3}
\begin{array}{ll}
\omega_2^1 = 0, & \omega_1^2 = 0,
\\
\omega_1^2 + \omega_1^1 -\omega_2^2 - \omega_2^1 = 0, &
 \omega_1^2 - \omega_1^1 + \omega_2^2 - \omega_2^1 = 0.
        \end{array}
\renewcommand{\arraystretch}{1}
$$

\textbf{5.} Finally we state the following theorem.

\begin{theorem}
Let $W (3, 2, 2)$ be a web with nonvanishing covector  $a \neq 0$ and
with integrable transversal $a$-distribution $\Delta$ $($conditions $(8)$ hold$)$.
All $2$-dimensional  webs $W (3, 2, 1)$ cut by the foliations of
a web $W (3, 2, 2)$  on the integral surfaces $V^2$
of  $\Delta$ are hexagonal if and only if the equations
\begin{equation}\label{eq:24}
 - b^i_{111} a_2^3 + 3  b^i_{(112)} a_2^2 a_1
-  3  b^i_{(122)} a_2 a_1^2 +  b^i_{222} a_1^3 = 0, \;\; i = 1, 2,
\end{equation}
hold.
\end{theorem}

\begin{corollary}
If the  transversal $a$-distribution  $\Delta$ of a web $W (3, 2, 2)$
 is the  distribution  defined by  the equations
 $\crazy{\omega}{\alpha}^1 = 0$, or $\crazy{\omega}{\alpha}^2 = 0$,
 or $\crazy{\omega}{\alpha}^1 + \crazy{\omega}{\alpha}^2 = 0$, or
 $\crazy{\omega}{\alpha}^1 - \crazy{\omega}{\alpha}^2 = 0$,
 then  all $2$-dimensional  webs $W (3, 2, 1)$ cut by the foliations of
a web $W (3, 2, 2)$  on the integral surfaces $V^2$
of  $\Delta$ are hexagonal if and only if respectively the following equations
\begin{equation}\label{eq:25}
 b^1_{222} = 0, \;\; b^2_{222} = 0,
\end{equation}

\begin{equation}\label{eq:26}
 b^1_{111} = 0, \;\; b^2_{111} = 0,
\end{equation}

\begin{equation}\label{eq:27}
\renewcommand{\arraystretch}{1.3}
\left\{
\begin{array}{ll}
- b^1_{111} + 3 (b^1_{(112)} - b^1_{(122)}) + b^1_{222} = 0,  \\
- b^2_{111} + 3 (b^2_{(112)} - b^2_{(122)}) + b^2_{222} = 0,
  \end{array}
\right.
\renewcommand{\arraystretch}{1}
 \end{equation}

\begin{equation}\label{eq:28}
\renewcommand{\arraystretch}{1.3}
\left\{
\begin{array}{ll}
 b^1_{111} + 3 (b^1_{(112)} + b^1_{(122)}) + b^1_{222} = 0,  \\
 b^2_{111} + 3 (b^2_{(112)} + b^2_{(122)}) + b^2_{222} = 0
  \end{array}
\right.
\renewcommand{\arraystretch}{1}
 \end{equation}
 hold.
\end{corollary}

\section{Nonisoclinic  webs $\boldsymbol{W} \boldsymbol{(} \boldsymbol{3},
\boldsymbol{2}, \boldsymbol{2}\boldsymbol{)}$}

\textbf{1.} For such a web, we write equations (1)--(6) and their prolongations:
\begin{equation}\label{eq:29}
-\crazy{\omega}{3}^i = \crazy{\omega}{1}^i + \crazy{\omega}{2}^i,
\end{equation}

\begin{equation}\label{eq:30}
\renewcommand{\arraystretch}{1.3}
\left\{
\begin{array}{ll}
  d \crazy{\omega}{1}^i = \crazy{\omega}{1}^j \wedge \omega_j^i +
         a_j \crazy{\omega}{1}^j \wedge \crazy{\omega}{1}^i, \\
         d\crazy{\omega}{2}^i = \crazy{\omega}{2}^j \wedge \omega_j^j -
         a_j \crazy{\omega}{2}^j \wedge \crazy{\omega}{2}^i,
        \end{array}
        \right.
\renewcommand{\arraystretch}{1}
\end{equation}

\begin{equation}\label{eq:31}
 d\omega_j^i - \omega_j^k \wedge \omega_k^i =
 b_{jkl}^i \crazy{\omega}{1}^k \wedge \crazy{\omega}{2}^l,
\end{equation}

\begin{equation}\label{eq:32}
 da_i - a_j \omega_i^j = p_{ij} \crazy{\omega}{1}^j
+  q_{ij}  \crazy{\omega}{2}^j,
\end{equation}

\begin{equation}\label{eq:33}
 b_{[j|l|k]}^i = \delta_{[k}^i p_{j]l} ,\;\;\;
 b_{[jk]l}^i = \delta_{[k}^i q_{j]l},
\end{equation}

\begin{equation}\label{eq:34}
 \nabla p_{ij} = \crazy{p}{1}_{ijk} \crazy{\omega}{1}^k +
 \crazy{p}{2}_{ijk} \crazy{\omega}{2}^k ,\;\;\;\;
 \nabla q_{ij} = \crazy{q}{1}_{ijk} \crazy{\omega}{1}^k +
 \crazy{q}{2}_{ijk} \crazy{\omega}{2}^k,
\end{equation}
where
$$
\nabla p_{ij} = dp_{ij} - p_{kj} \omega_i^k -
p_{ik}\omega_j^k ,\;\;\;\; \nabla q_{ij} = dq_{ij} - q_{kj}
\omega_i^k - q_{ik}\omega_j^k;
$$
\begin{equation}\label{eq:35}
\renewcommand{\arraystretch}{1.3}
\left\{
\begin{array}{ll}
          \crazy{p}{1}_{i[jk]}^{ } + p_{i[j}a_{k]} = 0, \;\;\;\;
          \crazy{q}{2}_{i[jk]}^{ } - q_{i[j}a_{k]} = 0 ,\\
          \crazy{p}{2}_{ijk}^{ } -  \crazy{q}{2}_{ijk}^{ } + a_m b_{ijk}^m = 0,
         \end{array}
         \right.
\renewcommand{\arraystretch}{1}
\end{equation}

\begin{equation}\label{eq:36}
 d p = p \omega_i^i + \crazy{p}{1}_i^{ } \crazy{\omega}{1}^i +
 \crazy{p}{2}_i^{ } \crazy{\omega}{2}^i ,\;\;\;
 d q = q \omega_i^i + \crazy{q}{1}_i^{ } \crazy{\omega}{1}^i +
 \crazy{q}{2}_i^{ } \crazy{\omega}{2}^i ,
\end{equation}
where
$$
 p_{[12]} = p \;\;\;\;\; q_{[12]} = q, \;\;\;\;
  \crazy{p}{k}_i^{ } = \crazy{p}{k}_{[12]i}^{ } ,\;\;\;\;
 \crazy{q}{k}_i^{ } = \crazy{q}{k}_{[12]i}^{ },
$$

\begin{equation}\label{eq:37}
 p^2 +  q^2 > 0,
\end{equation}

Recall that the condition (37) means that a web $W(3,2,2)$ is
nonisoclinic.

\textbf{2.} The next theorem gives analytic characterizations
for different types of  nonisoclinic webs $W (3, 2, 2)$.

\begin{theorem}
 For  nonisoclinic webs $W (3, 2, 2)$ of
different types we have the following analytic characterizations:
\begin{description}
\item[a)]
 A nonisoclinic web $W (3, 2, 2)$ is transversally geodesic if
and only if
\begin{equation}\label{eq:38}
 b_{(jkl)}^i = \delta^i_{(j} b_{kl)},
\end{equation}
where $b_{kl}$ is a $(0, 2)$-tensor.

\item[b)]
 A nonisoclinic web $W (3, 2, 2)$ is hexagonal
 if and only if
\begin{equation}\label{eq:39}
b^i_{(jkl)} = 0.
\end{equation}

\item[c)]
 A nonisoclinic web $W (3, 2, 2)$ is a
 Bol  web  $B_m$  if and only if
\begin{equation}\label{eq:40}
b^i_{j(kl)}  = 0.
\end{equation}

\item[d)]
A nonisoclinic web $W (3, 2, 2)$ is a group web if and only if
\begin{equation}\label{eq:41}
 b_{jkl}^i = 0,
\end{equation}

\item[e)]
A nonisoclinic three-web $W(3,2,2)$ can be uniquely extended
to a nonisoclinic web $AGW(4,2,2)$ $($see $[\mbox{{\rm G}}\; 88, 92])$
if and only if the following conditions hold:
\begin{equation}\label{eq:42}
 p \neq 0, \; q \neq 0, \;\; p \neq q,
\end{equation}
{\em and}
\begin{equation}\label{eq:43}
 \crazy{q}{ }_{ }^{ } (\crazy{q}{ }_{ }^{ }
 \crazy{p}{1}_i^{ } - \crazy{p}{ }_{ }^{ }
 \crazy{q}{1}_i^{ }) -\crazy{p}{ }_{ }^{ }
 (\crazy{q}{ }_{ }^{ }\crazy{p}{2}_i^{ } -
 \crazy{p}{ }_{ }^{ } \crazy{q}{2}_i^{ }) =
 \crazy{p}{ }_{ }^{ } \crazy{q}{ }^{ }_{ }
 (\crazy{p}{ }_{ }^{ } -\crazy{q}{ }_{ }^{ } )a_i.
\end{equation}
 The $4$th foliation of the web $AGW(4,2,2)$ is defined by
the equations
\begin{equation}\label{eq:44}
p \crazy{\omega}{1}^i + q \crazy{\omega}{2}^i = 0.
\end{equation}
\end{description}
\end{theorem}

Note that  part a) of Theorem 7 implies
 the following two very practical tests for a web
$W (3, 2, 2)$  to be nontransversally geodesic.

\begin{corollary}
\begin{description}
\item[a)]
The nonvanishing of any
of four components $b^i_{jjj}, \, i \neq j,$
implies that a web $W (3, 2, 2)$
is  not transversally geodesic.
\item[b)]
 If a web $W (3, 2, 2)$ is not a group $3$-web and
$b^1_{jkl} = 0$ or $b^2_{jkl} = 0$, then this web
is  not transversally geodesic.
\end{description}
\end{corollary}

{\sf Proof.} a) It follows from (38) that for a transversally geodesic
$b^i_{jjj} = 0, \; i \neq j$.

b) Suppose that  $b^2_{jkl} = 0$. Then it follows from (38)
that  $b_{jk} = 0$, and relation (38) implies that
 $b^1_{jkl} = 0$. As a result, the web $W (3, 2, 2)$
is a group 3-web. This contradicts to the conditions in b).
\rule{3mm}{3mm}

We describe another practical test  for a web
$W (3, 2, 2)$  to be transversally geodesic or not.
In general, to check whether a given web is transversally geodesic or
nontransversally geodesic, one may assume that the web is
transversally geodesic, calculate
$b_{ij}$ applying the formula $b_{ij} = \frac{3}{4} b^k_{(kij)}$
which is a consequence of (38), and substitute $b_{ij}$ obtained
into (38). If (38) will become the identity, the web
in question is transversally geodesic, and if (38) fails, the web
is nontransversally geodesic.

{\bf 3.} Suppose that 3 foliations $\lambda_1 ,\; \lambda_2$,
and $\lambda_3$ of a web $W (3, 2, 2)$  are given as the level
sets $u_{\xi}^i = \mbox{const.}\; (\xi = 1, 2, 3$) of the
following equations:
\begin{equation}\label{eq:45}
 \lambda_1 : u_1^i = x^i ;\;\;
 \lambda_2 : u_2^i = y^i ;\;\;
 \lambda_3 : u_3^i = f^i(x^j,y^k),\;\;\;\;
 i,j,k = 1,2.
\end{equation}

In order to characterize a web $W (3, 2, 2)$ given by
(45), we must find the \mbox{forms}
$\crazy{\omega}{\alpha}^i, \;\alpha = 1, 2, \;
\omega_j^i$, and the functions $a_{jk}^i ,\; b_{jkl}^i,\; a_i ,\;
p_{ij}, q_{ij}, \crazy{p}{1}_{ijk}, \crazy{p}{2}_{ijk},
\crazy{q}{1}_{ijk}, \linebreak \crazy{q}{2}_{ijk}, p, q, \crazy{p}{1}_{i},
\crazy{p}{2}_{i}, \crazy{q}{1}_{i}$, and $\crazy{q}{2}_{i}$.

The forms $\crazy{\omega}{\alpha}^i ,\; \omega_j^i$ and the
functions $a_{jk}^i$ and $b_{jkl}^i$ can be found my means of the
following formulas (see [AS 71], or [G 88], Section {\bf 8.1}, or
[AS 92], Section {\bf 1.6}):
\begin{equation}\label{eq:46}
 \crazy{\omega}{1}^i = \bar{f}_j^i dx^j ,\;\;\;
 \crazy{\omega}{2}^i = \tilde{f}_j^i dy^j ,\;\;\;
 \crazy{\omega}{3}^i = - du_3^i ,
\end{equation}
where $$
 \bar{f}_j^i = \displaystyle
\frac{\partial{f^i}}{\partial{x^j}}, \;\;\;\;
 \tilde{f}_j^i = \displaystyle
\frac{\partial{f^i}}{\partial{y^j}} , \;\;\;\;
 \det{(\bar{f}_j^i)} \neq 0 , \;\;\;
 \det{ (\tilde{f}_j^i)} \neq 0,
$$
and
\begin{equation}\label{eq:47}
 d\crazy{\omega}{1}^i = - d\crazy{\omega}{2}^i = \Gamma_{jk}^i
 \crazy{\omega}{1}^j \wedge \crazy{\omega}{2}^k ,
\end{equation}
\begin{equation}\label{eq:48}
 \Gamma_{jk}^i = - \displaystyle
\frac{\partial^2 f^i}{\partial{x^l} \partial{y^m}}
 \bar{g}_j^l \tilde{g}_k^m ,
\end{equation}
\begin{equation}\label{eq:49}
 \omega_j^i = \Gamma_{kj}^i \crazy{\omega}{1}^k + \Gamma_{jk}^i
 \crazy{\omega}{2}^k ,
\end{equation}
\begin{equation}\label{eq:50}
 a_{jk}^i = \Gamma_{[jk]}^i ,
\end{equation}
\begin{eqnarray}\label{eq:51}
\renewcommand{\arraystretch}{1.3}
 b_{jkl}^i & =  \displaystyle \frac{1}{2} \Biggl(
 \displaystyle
\frac{\partial{\Gamma_{kl}^i}}{\partial{x^m}} \bar{g}_j^m +
\displaystyle
 \frac{\partial{\Gamma_{jl}^i}}{\partial{x^m}} \bar{g}_k^m
- \displaystyle
 \frac{\partial{\Gamma_{kj}^i}}{\partial{y^m}} \tilde{g}_l^m
- \displaystyle \frac{\partial{\Gamma_{kl}^i}}{\partial{y^m}}
\tilde{g}_j^m \nonumber \\ &+   \Gamma_{jl}^m \Gamma_{km}^i -
\Gamma_{kj}^m \Gamma_{ml}^i
 + 2\Gamma_{kl}^m a_{mj}^i\Biggr) .
\renewcommand{\arraystretch}{1}
\end{eqnarray}
As to the functions $a_i ,\; p_{ij}, q_{ij},
\crazy{p}{1}_{ijk}, \crazy{p}{2}_{ijk},
\crazy{q}{1}_{ijk}, \crazy{q}{2}_{ijk}, p, q, \crazy{p}{1}_{i},
\crazy{p}{2}_{i}, \crazy{q}{1}_{i}$, and $\crazy{q}{2}_{i}$,
they can be easily calculated from equations
(30)--(36). Then we should check whether this web is
nonisoclinic, that is, whether condition (37) holds.
 Following this, we can investigate
to which of the classes indicated in Theorems 1, 3, 5, 7 the web in
question belongs. In the case e) of Theorem 7, we can find
equations (44) of
the 4th foliation of the web $W (4, 2, 2)$, an extension of the
3-web in question. Integrating these equations, we will find
closed-form equations of leaves of this 4th foliation.

In what follows, we will always assume that 3  foliations of a
web $W (3, 2, 2)$ are given as follows:
\begin{equation}\label{eq:52}
\renewcommand{\arraystretch}{1.3}
\left\{
        \begin{array}{ll}
         \lambda_1 : x^1 = \mbox{const.},\;\;\; x^2 =
\mbox{const.};\;\;\; \\
         \lambda_2 : y^1 = \mbox{const.},\;\;\; y^2 =
\mbox{const.}; \\ u^1_3 = f^1 (x^j, y^k) = \mbox{const.},\;\;\;
u_3^2 = f^2 (x^j, y^k) =\mbox{const.}
        \end{array}
        \right.
\renewcommand{\arraystretch}{1.3}
\end{equation}
In examples of webs, that we are going to present
in Section {\bf 3},  we will  only specify the functions $f^1
(x^j, y^k)$ and $f^2 (x^j, y^k)$.

Note that it follows from (49) that all forms $\omega_j^i$ are
expressed in terms of the forms $\crazy{\omega}{1}^i$ and
$\crazy{\omega}{2}^i$ only. This means that the forms $\pi_j^i =
\omega_j^i\Bigl|_{\crazy{\omega}{\alpha}^i = 0}$
 vanish, $\pi_j^i = 0$.
  It follows that the 2nd conditions in equations
  (9)--(12) are always valid for webs
  defined by equations (52), and the meaning of
  equations  (9)--(12) is that the transversal
  distribution $\Delta$ coincides with the distribution
  $\crazy{\omega}{1}^i = 0, \;\crazy{\omega}{2}^i = 0, \;
   \crazy{\omega}{1}^i  + \crazy{\omega}{2}^i = 0,
\;\crazy{\omega}{1}^i  - \crazy{\omega}{2}^i = 0$,
   respectively.

\textbf{4.} We will present a classification of nonisoclinic webs $W (3, 2, 2)$
given by equations (52).
\begin{description}
\item[{\bf A.}] Webs  with the
             integrable  transversal distribution $\Delta$.
     \begin{description} \item[${\bf A_{1}}.$] Webs  with the
             integrable  transversal distribution $a_1 \crazy{\omega}{1}^i
+ a_2 \crazy{\omega}{2}^i = 0, \; a_1, a_2 \neq 0 $ ((8) holds).
        \begin{description}
                                                      \item[${\bf A_{11}}$.]
                                             Webs for which the surfaces
                                           $V^2$ are geodesicly parallel ((17) holds).
 \item[{\bf A$_{12}$.}]
                                                     Webs   foliated into
                                                      2-dimensional  hexagonal webs
                                                     $W (3, 2, 1)$ ((8) and (24) hold).

                         \item[{\bf A$_{13}$.}] Webs with the integrable
                                  transversal distribution
    $\crazy{\omega}{\alpha}^1 + \crazy{\omega}{\alpha}^2   = 0$ ((11) and (15) hold).
         \begin{description}
\item[${\bf A_{131}}$.]
                                             Webs for which the surfaces
                                $V^2$ are geodesicly parallel ((11) and (22) hold).
 \item[${\bf A_{132}}$.]
                                                     Webs   foliated into
                                                      2-dimensional  hexagonal webs
                                              $W (3, 2, 1)$ ((11), (15) and (27) hold).
         \end{description}
   \item[${\bf A_{14}}$.] Webs with the integrable
                                  transversal distribution
    $\crazy{\omega}{\alpha}^1 - \crazy{\omega}{\alpha}^2   = 0$ ((12) and (16) hold).
         \begin{description}
\item[${\bf A_{141}}$.]
                                             Webs for which the surfaces
                                           $V^2$ are geodesicly parallel ((12)
                                           and (23) hold).
 \item[${\bf A_{142}}$.]
                                                     Webs   foliated into
                                                      2-dimensional  hexagonal webs
                                                     $W (3, 2, 1)$ ((12), (16), and (28) hold).
         \end{description}
\end{description}
           \item[${\bf A_{2}}.$] Webs  with the
            integrable  transversal distribution $\crazy{\omega}{\alpha}^1 = 0$
                  ((9) and (13) hold).
  \begin{description}
\item[${\bf A_{21}}$.]
                                             Webs for which the surfaces
                                      $V^2$ are geodesicly parallel ((9) and (20) hold).
 \item[${\bf A_{22}}$.]
                                       Webs   foliated into
                                                      2-dimensional  hexagonal webs
                                                     $W (3, 2, 1)$ ((9), (13), and (25) hold).
     \end{description}
            \item[${\bf A_3}$.] Webs  with the integrable  transversal
                                     distribution $\crazy{\omega}{\alpha}^2 = 0$
                                     ((10) and (14) hold).
   \begin{description}
                               \item[${\bf A_{31}}$.]
                                             Webs for which the surfaces
                                           $V^2$ are geodesicly parallel ((10) and (21) hold).
 \item[${\bf A_{32}}$.]
                                Webs   foliated into
                              2-dimensional  hexagonal webs
                             $W (3, 2, 1)$ ((10), (14), and (26) hold).
  \end{description}
   \end{description}
        \item[{\bf B.}] Webs with nonintegrable  transversal
                          distribution $\Delta$.
\item[{\bf C.}] Nontransversally geodesic webs ((38) does not hold).
\item[{\bf D.}] Transversally geodesic webs ((38) holds).
\begin{description}
\item[{\bf D${}_1$.}] Hexagonal webs ((39) holds).
\begin{description}
\item[{\bf D${}_{11}$.}] Bol webs ((40) holds).
\item[{\bf D${}_{12}$.}] Group webs ((41) holds).
\end{description}
\end{description}
\item[{\bf E.}] Webs with different relations among $p_{ij}$ and
                   $q_{ij}$.

                   For webs defined by equa\-tions
                   (52), each component of the tensors $p_{ij}$ and
                   $q_{ij}$ is an absolute invariant, and the vanishing
                   any of these components distinguishes a class of 3-webs
                   $W (3, 2, 2)$.  Note  that since we consider only
                   nonisoclinic webs $W (3, 2, 2)$, the tensors $p_{ij}$ and
                   $q_{ij}$ cannot be simultaneously symmetric.
                   Note also that the conditions $a_2 = 0, \;\;
                   \omega_2^1 = 0$ imply the conditions $p_{2i} =
                   q_{2i} = 0$ (cf. (20)), and  the conditions $a_1 = 0, \;\;
                   \omega_1^2 = 0$ imply the conditions $p_{1i} =
                   q_{1i} = 0$ (cf. (21)).
                   We indicate some of classes of such webs.

\begin{description}
\item[${\bf E_1}$.] Webs with $ a_1 \neq 0, \; a_2  \neq 0, \; a_1 \neq a_2$.
\begin{description}
\item[${\bf E_{11}}$.] Webs with $p_{i2} = q_{i1} = 0$.
\begin{description}
\item[${\bf E_{111}}$.] Webs with  $p_{21} =  - q_{12}$.
\end{description}
\item[${\bf E_{12}}$.] Webs with $p_{12} = p_{21}, \;
q_{11} = q_{22}, \; q_{12} = -q_{21}$.
\item[${\bf E_{13}}$.] Webs with $p_{11} = p_{22} = 0$.
\begin{description}
\item[${\bf E_{131}}$.] Webs with $p_{12} = q_{12}, \;
p_{21} = q_{21}$.
\end{description}
\end{description}
\item[${\bf E_2}$.] Webs with
%$a_2 = 0, \; \omega_2^1 = 0$ (or
$p_{2i} = q_{2i} = 0$.
\begin{description}
\item[${\bf E_{21}}$.] Webs with $q_{12} = 0$.
\item[${\bf E_{22}}$.] Webs with $p_{1i} = 0$.
\item[${\bf E_{23}}$.] Webs with $p_{12} = q_{12}$.
\end{description}
\item[${\bf E_3}$.] Webs with
%$a_1 = 0, \; \omega_1^2 = 0$ (or
$p_{1i} = q_{1i} = 0$.
\begin{description}
\item[${\bf E_{31}}$.] Webs with $q_{2j} = 0$.
\item[${\bf E_{32}}$.] Webs with $p_{21} = q_{21}$.
\begin{description}
\item[${\bf E_{321}}$.] Webs with $p_{22} = 0$.
\end{description}
\item[${\bf E_{33}}$.] Webs with $p_{21} = -q_{21}$.
\end{description}
\end{description}

\item[{\bf F.}] Webs extendable to exceptional  webs $W (4, 2, 2)$ of maximum 2-rank
((42) and (43) hold).
\item[{\bf G.}] Webs nonextendable to exceptional webs $W (4, 2, 2)$ of
 maximum 2-rank ((42) or (43) does not hold).
 \begin{description}
\item[${\bf G{}_{1}}$.] Webs with $p = 0, \; q \neq 0$.
\item[${\bf G{}_{2}}$.] Webs with $q = 0, \; p \neq 0$.
\item[${\bf G{}_{3}}$.] Webs with $p = q \neq 0$.
\item[${\bf G{}_{4}}$.] Webs with $p \neq 0, \; q \neq 0, \; p \neq q$
and for which condition (43) does not hold.
\end{description}
\end{description}

\rem{The classification presented above is complete in the sense
that any web $W (3, 2, 2)$ belongs to the class {\bf A} or  {\bf B},
 {\bf C} or  {\bf D}, {\bf F} or  {\bf G},  {\bf E} or  the class of
 webs with other or no connections among $p_{ij}$ and $q_{ij}$.
}

\rem{
It is easy to see a geometric meaning of the  classes
${\bf G_1-G_3}$ for
which conditions (42) do not hold. By (44), the 4th foliation
$\lambda_4$ defining a nonisoclinic web $W (4, 2, 2)$ coincides
with  the foliations $\lambda_1, \;  \lambda_2$, and $\lambda_3$,
respectively.
As to the  class ${\bf G_4}$, the 4th foliation is well defined but not
integrable.
}

\rem{It is easy to give  characterizations of  webs of the classes
${\bf  E_{1}-  E_{3}}$:

\begin{itemize}
\item Class ${\bf  E_{11}}$: for webs of this class, both components,
 $a_1$ and $a_2$, are covariantly
      constant on the 2-dimensional distribution defined by the  equations
      $\crazy{\omega}{1}^1 = \crazy{\omega}{2}^2 =  0$.

\begin{itemize}
\item Class ${\bf  E_{111}}$: for webs of this class, in addition
to be covariantly  constant on the 2-dimensional distribution defined by the equations
      $\crazy{\omega}{1}^1 = \crazy{\omega}{2}^2  = 0$,  the components $a_1$ and
$a_2$  satisfy
      the exterior quadratic equation $ \nabla a_1 \wedge \crazy{\omega}{1}^1
      = \nabla a_2 \wedge \crazy{\omega}{2}^2$.
On some instances (see, for example Example {\bf 11} below) this
condition is necessary and sufficient for integrability of the
distribution defined by the equation $a_1 \crazy{\omega}{1}^1
= a_2 \crazy{\omega}{2}^2$.

\end{itemize}

\item Class ${\bf  E_{12}}$:  for webs of this class, the quantities $a_1$
      and $a_2$ satisfy the following 3 exterior quartic equations:
      $$
\renewcommand{\arraystretch}{1.3}
\left\{
\begin{array}{ll}
\nabla a_1 \wedge \crazy{\omega}{1}^1 \wedge \crazy{\omega}{2}^1 \wedge \crazy{\omega}{2}^2 =
-\nabla a_2 \wedge \crazy{\omega}{1}^2 \wedge \crazy{\omega}{2}^1 \wedge \crazy{\omega}{2}^2, \\
\nabla a_1 \wedge \crazy{\omega}{1}^1 \wedge \crazy{\omega}{1}^2 \wedge \crazy{\omega}{2}^2 =
-\nabla a_2 \wedge \crazy{\omega}{1}^1 \wedge \crazy{\omega}{1}^2 \wedge \crazy{\omega}{2}^1, \\
\nabla a_1 \wedge \crazy{\omega}{1}^1 \wedge \crazy{\omega}{1}^2 \wedge \crazy{\omega}{2}^1 =
\nabla a_2 \wedge \crazy{\omega}{1}^1 \wedge \crazy{\omega}{1}^2 \wedge \crazy{\omega}{2}^2.
        \end{array}
        \right.
\renewcommand{\arraystretch}{1}
$$

\item Class ${\bf  E_{13}}$: for webs of this class, the component $a_1$
is covariantly constant on the one-dimensional distribution defined by the equation
      $\crazy{\omega}{1}^2 =  \crazy{\omega}{2}^i =0$, and the component $a_2$
is covariantly constant on the 2-dimensional distribution defined by the
equations   $\crazy{\omega}{2}^1 =  \crazy{\omega}{2}^i = 0$.
\begin{itemize}
\item Class ${\bf  E_{131}}$: for webs of this class, both the components $a_1$
and $a_2$ are  covariantly constant on the 2-dimensional distributions defined
by the equations $\crazy{\omega}{2}^1 =  \crazy{\omega}{1}^2 =\crazy{\omega}{2}^2 =0$
      and  $\crazy{\omega}{1}^1 =  \crazy{\omega}{1}^1
      =\crazy{\omega}{2}^2 =0$, respectively.
\end{itemize}

\item Class ${\bf  E_{21}}$: for webs of this class,  the component $a_1$ is covariantly
      constant on the ome-dimensional distribution defined by the equations
      $\crazy{\omega}{1}^i = \crazy{\omega}{2}^1 = 0$.

\item Class ${\bf  E_{22}}$: for webs of this class,  the component $a_1$ is covariantly
      constant on the 2-dimensional distribution defined by the equations
      $\crazy{\omega}{2}^i = 0$.

\item Class ${\bf  E_{23}}$: for webs of this class, the component $a_1$ is covariantly
      constant on the one-dimensional distribution defined by the equations
      $\crazy{\omega}{1}^1 = \crazy{\omega}{2}^1
      = \crazy{\omega}{1}^2  + \crazy{\omega}{2}^2 = 0$.

\item Class ${\bf  E_{31}}$: for webs of this class, the component $a_1$ is covariantly
      constant on  the entire web, and the component $a_2$ is covariantly
      constant on the 2-dimensional distribution defined by the equations
      $\crazy{\omega}{1}^i = 0$.

\item Class ${\bf  E_{32}}$: for webs of this class, the component $a_1$ is covariantly
      constant on the entire web, and the component $a_2$ is covariantly
      constant on  the one-dimensional distribution defined by the equations
      $\crazy{\omega}{1}^2 = \crazy{\omega}{2}^2 = \crazy{\omega}{1}^1  + \crazy{\omega}{2}^1 = 0$.

\begin{itemize}
 \item Class ${\bf  E_{321}}$: for webs of this class, the component $a_2$ is covariantly
      constant on  the 2-dimensional distribution defined by the equations
      $ \crazy{\omega}{2}^2 = \crazy{\omega}{1}^1  + \crazy{\omega}{2}^1 = 0$.
\end{itemize}
\end{itemize}
}

\rem{
If a web is given by equation (52), then
the webs of the classes ${\bf A_1, A_2,
A_3, A_{13},  A_{14}}$ as well as webs of the classes
${\bf A_{11},  A_{131},  A_{141},  A_{21},
 A_{31}}$ and the classes ${\bf A_{12},  A_{132},
 A_{142},  A_{22},  A_{32}}$
 could be equivalent one to another.
They are different by a location of the integrable transversal distribution
$\Delta$. In fact, if there is no additional conditions on webs, then
by transformations $x^i = \phi^i (x^j), \; y^i = \psi^i (y^j)$, we
can make the integrable transversal distributions $\Delta$ of any two
of them to coincide. However, in our examples in Section {\bf 3},
we always have $\pi_i^j = 0$, and  above mentioned
specializations are impossible.
In addition, in our examples,
there will be additional conditions on webs, and in general, under
the above mentioned transformation (if it would be possible),
these additional conditions for the first transformed web
 do not coincide with the additional conditions for the 2nd
 web. These are the reasons that in our classification, we consider
all above mentioned classes as different ones.
}
\rem{The classification presented above has some overlapping classes:
for example, the classes ${\bf A_{131}}$ and  ${\bf A_{141}}$ are
 subclasses of the class  ${\bf A_{11}}$,
 the classes ${\bf A_{132}}$ and  ${\bf A_{142}}$ are
 subclasses of the class  ${\bf A_{12}}$, and
 the classes  ${\bf A_{21}}$ and   ${\bf A_{31}}$
 are subclasses of the classes ${\bf E_{2}}$ and ${\bf E_{3}}$,
 respectively.
}

\rem{In general, we must prove the existence theorem for all
the classes listed above. Such a theorem can be proved
for a web of general kind of each class using the well-known
Cartan's test or it can be proved by finding examples
of webs of these classes. The examples
of webs in Section {\bf 3}  prove the existence of webs of most of the classes
of our classification.
}

\section{Examples of nonisoclinic  webs $W (3, 2, 2)$}

\setcounter{theorem}{0}

In subsections {\bf 1--5}, we will consider examples
of nonisoclinic webs $W (3, 2, 2)$ of the classes
${\bf F}$ and ${\bf G_1-G_4}$,
and for each of them, we will indicate to which other  classes
it belongs.

{\bf 1.} {\bf An example of an extendable nonisoclinic web $W (3,
2, 2)$ (Class {\bf F})}. We start from an example of a nonisoclinic web $W (3, 2,
2)$ that can be extended to a nonisoclinic web $\mathcal{G}_3 (4, 2, 2)$.

\examp{\label{examp:1} When the author was constructing the
exceptional  web $\mathcal{G}_3 (4, 2, 2)$ (see [G 87]; see also [G 88],
Example {\bf 8.4.5}, p. 431; [AG 96], Example {\bf 5.5.3}, p. 201;
[AG 00], Section {\bf 1.4}, p. 101; and [AS 92], Ch. 8, problem {\bf 6},
p. 307), he started from an example of
a web $W (3, 2, 2)$  defined  by the equations
\begin{equation}\label{eq:53}
u_3^1 = x^1 + y^1 + \displaystyle \frac{1}{2} (x^1)^2 y^2, \;\;
 u_3^2 = x^2 + y^2 - \displaystyle \frac{1}{2} x^1 ( y^2)^2
\end{equation}
in a domain $x^1 y^2 \neq \pm 1$,  where $\Delta = 1 + x^1 y^2$.

Using (48)--(51) and (32)--(37), we find that
\begin{equation}\label{eq:54}
\Gamma^1_{12} = - \displaystyle \frac{x^1}{\Delta (2 - \Delta)},
\;\;
\Gamma^2_{12} =  \displaystyle \frac{y^2}{\Delta (2 -\Delta)},
\;\;
\Gamma^1_{11} = \Gamma^1_{2i} = \Gamma^2_{22} = \Gamma^2_{i1} = 0,
\end{equation}
\begin{equation}\label{eq:55}
\renewcommand{\arraystretch}{1.3}
\left\{
\begin{array}{ll}
\omega^1_1 = -\displaystyle \frac{x^1}{\Delta (2 - \Delta)} \crazy{\omega}{2}^2, &
\omega^1_2 =  -\displaystyle \frac{x^1}{\Delta (2 - \Delta)}
\crazy{\omega}{1}^1,\\
\omega^2_1 = \displaystyle \frac{y^2}{\Delta (2 -\Delta)} \crazy{\omega}{2}^2, &
\omega^2_2 = \displaystyle \frac{y^2}{\Delta (2 -\Delta)} \crazy{\omega}{1}^1,
        \end{array}
        \right.
\renewcommand{\arraystretch}{1}
\end{equation}
 \begin{equation}\label{eq:56}
a_1 = \displaystyle \frac{y^2}{\Delta (2 - \Delta)},\;\;\;\;
a_2 = \displaystyle \frac{x^1}{\Delta (2 - \Delta)},
\end{equation}
\begin{equation}\label{eq:57}
\renewcommand{\arraystretch}{1.3}
\left\{
\begin{array}{ll}
p_{11}=\displaystyle \frac{2(y^2)^2 (\Delta -1)}{\Delta^3 (2 -
\Delta)^2}, &
    p_{21}=\displaystyle
\frac{\Delta^2 - 2\Delta + 2}{\Delta^3 (2 - \Delta)^2},\\
    q_{22}= \displaystyle
\frac{2(x^1)^2 (\Delta -1)}{\Delta^2 (2 - \Delta)^3}, &
   q_{12}=\displaystyle
\frac{\Delta^2 - 2\Delta + 2}{\Delta^2 (2 - \Delta)^3}, \\
         p_{i2} = q_{i1}= 0, &
        \end{array}
        \right.
\renewcommand{\arraystretch}{1}
\end{equation}
\begin{equation}\label{eq:58}
 p = - \displaystyle
\frac{\Delta^2 - 2\Delta + 2}{2\Delta^3 (2 - \Delta)^2}, \;\;\;\;
 q =  \displaystyle
\frac{\Delta^2 - 2\Delta + 2}{2\Delta^2 (2 - \Delta)^3},
\end{equation}
\begin{equation}\label{eq:59}
\renewcommand{\arraystretch}{1.3}
\left\{
\begin{array}{ll}
         \crazy{p}{1}_1 = \displaystyle
\frac{y^2(-\Delta^3 + 4\Delta^2 - 8\Delta +6)}{\Delta^5
         (2 - \Delta)^3}, &
         \crazy{p}{2}_2 = \displaystyle
\frac{x^1(-\Delta^3 + 3\Delta^2 - 6\Delta +4)}{\Delta^4
         (2 - \Delta)^4}, \\
         \crazy{q}{1}_1 = \displaystyle
\frac{y^2(\Delta^3 - 3\Delta^2 + 6\Delta -4)}{\Delta^4
         (2 - \Delta)^4}, &
         \crazy{q}{2}_2 = \displaystyle
\frac{x^1(\Delta^3 - 2\Delta^2 + 4\Delta -2)}{\Delta^3
         (2 - \Delta)^5}, \\
    \crazy{p}{1}_2 = \crazy{q}{1}_2 = \crazy{p}{2}_1
= \crazy{q}{2}_1 = 0, &
        \end{array}
        \right.
\renewcommand{\arraystretch}{1.3}
\end{equation}
\begin{equation}\label{eq:60}
\renewcommand{\arraystretch}{1.3}
\left\{
\begin{array}{ll}
b^1_{112} = -\displaystyle \frac{\Delta^2 - 2 \Delta + 2}{\Delta^3 (2 - \Delta)^2}, &
b^1_{111} = b^1_{12i} =  0,\\
b^1_{212} = -  \displaystyle \frac{2(x^1)^2 (\Delta - 1)}{\Delta^2 (2 - \Delta)^3}, &
b^1_{211} = b^1_{22i} = 0,\\
b^2_{112} = -  \displaystyle \frac{2(y^2)^2 (\Delta - 1)}{\Delta^3 (2 -\Delta)^2}, &
b^2_{111} = b^2_{121} = b^2_{122} = 0,\\
b^2_{212} = - \displaystyle \frac{\Delta^2 - 2 \Delta + 2}{\Delta^2 (2 - \Delta)^3}, &
b^2_{211} = b^2_{22i} = 0.
        \end{array}
        \right.
\renewcommand{\arraystretch}{1}
\end{equation}

It follows from  (56) and (57) that the conditions (8) do not hold.
Therefore, the $a$-distribution defined by equations (7)
is not integrable, and the web (53) belongs to the class {\bf B}.

It follows from (60) and (38) that  the web
(53) is not transversally geodesic (and consequently is neither
hexagonal, nor Bol, nor group web). In fact, if (38) holds, then
(60) gives $b^1_{111} = b_{11} = 0$ but it follows from (38) that
$b_{11} = - \frac{1}{12} b^2_{112} \neq 0$. This contradiction proves that
the web (53) belongs to the class {\bf C}.
 By (57), the web (53) belongs to the class ${\bf E_{11}}$.

 It follows from  (56), (58), and (59) that
the conditions (42) and (43) are satisfied. Hence the equations
(53) define a {\it nonisoclinic $3$-web $W(3,2,2)$ that
can be extended to a nonisoclinic $4$-web $\mathcal{G}_3 (4,2,2)$.} Thus
the web (53) belongs to the class {\bf F}.

As we proved in [AG 87] (see also [G 88], Section {\bf 8.4}), the
leaves of the 4th foliation of the web $W (4, 2, 2)$ are
defined as level sets of the following functions:
\begin{equation}\label{eq:61}
 u_4^1 = -x^1 + y^1 + \displaystyle
\frac{1}{2} (x^1)^2 y^2,\;\;\;\;
 u_4^2 = x^2 -y^2 -\displaystyle
\frac{1}{2} x^1 (y^2)^2.
\end{equation}
It is also shown in [AG 87] (see also [G 88], Section {\bf 8.4})
that  the 4-web defined by (53) and (61) is of maximum $2$-rank, and
the only abelian 2-equation has the form
\begin{equation}\label{eq:62}
 2dx^1 \wedge dx^2 + 2dy^1 \wedge dy^2 - du_3^1 \wedge du_3^2 +
 du_4^1 \wedge du_4^2 = 0.
\end{equation}

Note that {\it the $4$-web constructed in this example represents
the first and only known example of a nonisoclinic web $W(4,2,2)$
of maximum {\rm 2}-rank.}
}

\textbf{2. Examples of nonextendable
 nonisoclinic  webs} $\boldsymbol{W}
 (\boldsymbol{3}, \boldsymbol{2}, \boldsymbol{2})$
 \textbf{with}
 $\boldsymbol{p} \boldsymbol{=} \boldsymbol{0}, \;
 \boldsymbol{q} \boldsymbol{\neq} \boldsymbol{0}$
\textbf{(Class}
 $\bf{G_1).}$

\examp{\label{examp:2} Suppose that a web $W(3,2,2)$ is
given by
\begin{equation}\label{eq:63}
 u_3^1 = x^1 + y^1,\;\;\;\;
 u_3^2 = x^1 y^1 + x^2 y^2
\end{equation}
in a domain $x^2 \neq 0, \; y^2 \neq 0$
(see [G 88], Example {\bf 8.1.28}, p. 392).
In this case using (48)--(51) and (32)--(37), we find that
\begin{equation}\label{eq:64}
\Gamma^2_{11} = \displaystyle \frac{y^1}{y^2} - 1, \;\;
\Gamma^2_{21} =   \displaystyle  \frac{1}{y^2}, \;\;
\Gamma^2_{i2} =  \Gamma^1_{jk} = 0,
\end{equation}
\begin{equation}\label{eq:65}
\omega^2_1 = \Bigl(\frac{y^1}{y^2}-1\Bigr) (\crazy{\omega}{1}^1 +
\crazy{\omega}{2}^1) + \frac{1}{y^2} \crazy{\omega}{1}^2 ,\;\;
\omega^2_2 =   \frac{1}{y^2} \crazy{\omega}{2}^1, \;\;
\omega^1_i = 0,
\end{equation}
\begin{equation}\label{eq:66}
 a_1 = -\displaystyle \frac{1}{y^2},\;\;\;  a_2 = 0,
\end{equation}
\begin{equation}\label{eq:67}
  p_{ij} = q_{2i}  = 0, \;\;
q_{11} = -\displaystyle \frac{x^1}{x^2 (y^2)^2}, \;\;
q_{12} = -\displaystyle \frac{1}{x^2 (y^2)^2},
\end{equation}

\vspace*{-5mm}

\begin{equation}\label{eq:68}
 p = 0, \;\;\; q = -\displaystyle
\frac{1}{2 x^2 (y^2)^2},
\end{equation}
\begin{equation}\label{eq:69}
\renewcommand{\arraystretch}{1.3}
\left\{
\begin{array}{ll}
b^1_{jkl} = 0, \;\;
b^2_{112} = \displaystyle \frac{y^1}{x^2 (y^2)^2}, \;\;
b^2_{121} = -\displaystyle \frac{x^1}{x^2 (y^2)^2}, \;\;
b^2_{211} = \frac{y^1 - x^1}{x^2 (y^2)^2}, \\
b^2_{111} = -\displaystyle \frac{1}{y^2} -\displaystyle  \frac{x^1 y^1}{x^2 (y^2)^2},\;\;
b^2_{122} = b^2_{221} = \displaystyle \frac{1}{x^2(y^2)^2}, \;\;
b^2_{2i2} = 0.
        \end{array}
        \right.
\renewcommand{\arraystretch}{1}
\end{equation}

Equations (67) and (68) show that the web (63) belongs to the
classes ${\bf E_{22}}$ and ${\bf G_{1}}$.
Since by (69), $b^2_{111} \neq 0$, it follows from  (38) that
the web (63) is not transversally geodesic
(and consequently is neither  hexagonal, nor Bol, nor group web).
Thus it belongs to the class  {\bf C}.
Equations (66), (67), (69), (9), (13), (20), and (25)
 show that
this web belongs to the classes ${\bf A_{21}}$
 and ${\bf A_{22}}$.
}

\examp{\label{examp:3} Consider the web $W(3,2,2)$ given by
\begin{equation}\label{eq:70}
 u_3^1 = x^1 y^1 - x^2 y^2 ,\;\;\;
 u_3^2 = x^1 y^2 + x^2 y^1
\end{equation}
in a domain where $y^1$ and $y^2$ are not 0 simultaneously,
and $x^1$ and $x^2$ are not 0 simultaneously
(see [G 88], Example {\bf 8.1.26}, p. 391 or [G 92], p. 341).
In this case using (48)--(51) and (32)--(37), we find that
\begin{equation}\label{eq:71}
\renewcommand{\arraystretch}{1.3}
\left\{
\begin{array}{ll}
\Gamma^1_{11}  = \Gamma^1_{22} = - \Gamma^2_{12} =  \Gamma^2_{21}
   = - \displaystyle \frac{\alpha}{\Delta_1 \Delta_2}, \\
\Gamma^1_{12} = -\Gamma^1_{21} = \Gamma^2_{11} = \Gamma^2_{22}
         = \displaystyle \frac{\beta}{\Delta_1 \Delta_2},
  \end{array}
   \right.
\renewcommand{\arraystretch}{1}
\end{equation}
\begin{equation}\label{eq:72}
\renewcommand{\arraystretch}{1.3}
\left\{
\begin{array}{ll}
\omega^1_1 =  \displaystyle \frac{1}{\Delta_1 \Delta_2}
\bigl[ - \alpha(\crazy{\omega}{1}^1 + \crazy{\omega}{2}^1)
+ \beta(\crazy{\omega}{2}^2 - \crazy{\omega}{1}^2) \bigr], \\
\omega^1_2 =  \displaystyle \frac{1}{\Delta_1 \Delta_2}
\bigl[\beta(\crazy{\omega}{1}^1 - \crazy{\omega}{2}^1)
- \alpha(\crazy{\omega}{1}^2 + \crazy{\omega}{2}^2) \bigr], \\
\omega^2_1 =  \displaystyle \frac{1}{\Delta_1 \Delta_2}
\bigl[\beta(\crazy{\omega}{1}^1 + \crazy{\omega}{2}^1)
+ \alpha(\crazy{\omega}{2}^2 - \crazy{\omega}{1}^2) \bigr], \\
\omega^2_2 =  \displaystyle \frac{1}{\Delta_1 \Delta_2}
\bigl[\alpha (\crazy{\omega}{1}^1 - \crazy{\omega}{2}^1)
+ \beta(\crazy{\omega}{1}^2 + \crazy{\omega}{2}^2)\bigr],
        \end{array}
       \right.
\renewcommand{\arraystretch}{1}
\end{equation}
\begin{equation}\label{eq:73}
 a_1  =  \displaystyle
\frac{2\alpha}{\Delta_1 \Delta_2},  \;\;\; a_2  =
\displaystyle -\frac{2\beta}{\Delta_1 \Delta_2},
\end{equation}
\begin{equation}\label{eq:74}
\renewcommand{\arraystretch}{1.3}
\left\{
\begin{array}{ll}
p_{11} = \displaystyle \frac{4 (y^2)^2}{\Delta_1 \Delta_2^2}, \;\;
p_{22} = \displaystyle  \frac{4 (x^1)^2}{\Delta_1 \Delta_2^2},&
q_{11} =  q_{22} = \displaystyle \frac{4 (y^2)^2}{\Delta_1^2 \Delta_2}, \\
p_{12} = p_{21} = - \displaystyle  \frac{4 x^1 x^2}{\Delta_1 \Delta_2^2}, &
q_{12} = - q_{21} = - \displaystyle \frac{4y^1 y^2}{\Delta_1^2 \Delta_2},
\end{array}
\right.
\renewcommand{\arraystretch}{1}
\end{equation}
\begin{equation}\label{eq:75}
p = 0, \;\;\; q = - \displaystyle \frac{4y^1 y^2}{\Delta_1^2 \Delta_2},
\end{equation}
\begin{equation}\label{eq:76}
\renewcommand{\arraystretch}{1.3}
\left\{
\begin{array}{llll}
b^1_{111} = \displaystyle \frac{2 \beta u_3^2}{\Delta_1^2 \Delta_2^2},&
b^1_{112} = \displaystyle \frac{2x^1 y^1}{\Delta_1 \Delta_2^2}, &
b^1_{121} = -\displaystyle \frac{2y^1 y^2}{\Delta_1^2 \Delta_2}, &
b^1_{211} = 0, \\
b^1_{222} = -\displaystyle \frac{2 \beta u_3^1}{\Delta_1^2 \Delta_2^2},&
b^1_{122} = - \displaystyle \frac{1}{\Delta_1 \Delta_2}, &
b^1_{212} = \displaystyle \frac{2(y^2)^2}{\Delta_1^2 \Delta_2}, &
b^1_{221} = \displaystyle \frac{(x^1)^2 - (x^2)^2}{\Delta_1 \Delta_2^2},\\
b^2_{111} = \displaystyle \frac{2 \alpha u_3^2}{\Delta_1^2 \Delta_2^2},&
b^2_{112} = \displaystyle \frac{2(x^2)^2}{\Delta_1 \Delta_2^2},&
b^2_{121} = \displaystyle \frac{2x^2 \gamma}{\Delta_1^2 \Delta_2^2},&
b^2_{211} = \displaystyle \frac{2x^1 \rho}{\Delta_1^2 \Delta_2^2}, \\
b^2_{222} = -\displaystyle \frac{2 \alpha u_3^1}{\Delta_1^2 \Delta_2^2},&
b^2_{122} = \displaystyle \frac{\alpha^2}{\Delta_1^2 \Delta_2^2},&
b^2_{212}= - \displaystyle \frac{2x^2 \gamma}{\Delta_1^2 \Delta_2^2},&
 b^2_{221}= \displaystyle \frac{2y^2 \sigma}{\Delta_1^2 \Delta_2^2}.
        \end{array}
        \right.
\renewcommand{\arraystretch}{1}
\end{equation}
where
$$
\renewcommand{\arraystretch}{1.3}
\left\{
\begin{array}{ll}
\Delta_1 = (y^1)^2 + (y^2)^2 , & \Delta_2 = (x^1)^2 + (x^2)^2, \\
 \alpha = x^1 y^1 + x^2 y^2, & \beta = x^1 y^2 - x^2 y^1,\\
 \gamma = x^2 [(y^1)^2 - (y^2)^2] + 2 x^1 y^1 y^2, &
 \rho = x^1 [(y^1)^2 - (y^2)^2] + 2 x^2 y^1 y^2, \\
 \sigma  = y^1 [(x^1)^2 - (x^2)^2] + 2 x^1 x^2 y^2. &
\end{array}
\right.
\renewcommand{\arraystretch}{1}
$$

It follows from  (73) and (74) that

$$
a_2^2 q_{11} -  2a_1 a_2 q_{(12)} + a_1^2 q_{22}
= \displaystyle \frac{16 \beta^2 (y^2)^2}{\Delta_1^3 \Delta^4_2} \neq 0,
$$
i.e,  the conditions (8) do not hold.
Therefore, the $a$-distribution  defined by equations
(7) is not integrable, and the web (70) belongs to the class {\bf B}.

Since by (76), $b^1_{222} \neq 0$,   the web
(70) is not transversally geodesic.
 (and consequently is neither
 hexagonal, nor Bol,  nor group web).
So this web
belongs to the class {\bf C}.  By (74) and (75), the web (70)
belongs to the classes  ${\bf E_{12}}$ and ${\bf G_1}$.
}

\textbf{ 3.    Example of a nonextendable nonisoclinic
web} $\boldsymbol{W} \boldsymbol{(}\boldsymbol{3},
\boldsymbol{2}, \boldsymbol{2}\boldsymbol{)}$
\textbf{with}  $\boldsymbol{p} \boldsymbol{\neq}
\boldsymbol{0}, \; \boldsymbol{q} \boldsymbol{=} \boldsymbol{0}$
\textbf{(Class} ${\bf G_2)}$.

\examp{\label{examp:4} A web $W(3,2,2)$ is given by the
equations
\begin{equation}\label{eq:77}
 u_3^1 = x^1 + y^1,\;\;\;
 u_3^2 = (x^1)^2 y^2 + (x^2)^2 y^1
\end{equation}
in a domain of ${\bf R}^4$ where $x^1, x^2, y^1 \neq 0$.

In this case using (48)--(51) and (32)--(37), we find that
\begin{equation}\label{eq:78}
\Gamma^1_{jk} = \Gamma^2_{22} = 0, \;\;
\Gamma^2_{11} =
\displaystyle \frac{(y^2)^2}{(x^1)^2} +
\displaystyle \frac{x^1 y^2}{x^2 y^1}, \;\;
\Gamma^2_{12} = -\displaystyle \frac{1}{(x^1)^2}, \;\;
\Gamma^2_{21} = - \displaystyle \frac{1}{2 x^2 y^1},
\end{equation}
\begin{equation}\label{eq:79}
\renewcommand{\arraystretch}{1.3}
\left\{
\begin{array}{ll}
\omega_i^1 =  0, \;\;
\omega^2_2 = - \displaystyle \frac{1}{(x^1)^2}\crazy{\omega}{1}^1
- \displaystyle \frac{1}{2x^2 y^1}\crazy{\omega}{2}^1, \\
\omega^2_1 =  \Bigl(\displaystyle \frac{(x^2)^2}{(x^1)^2}
           + \displaystyle \frac{x^1 y^2}{x^2   y^1}\Bigr)
           (\crazy{\omega}{1}^1 + \crazy{\omega}{2}^1)
           -  \displaystyle \frac{1}{2x^2   y^1} \crazy{\omega}{1}^2
           -  \displaystyle \frac{1}{(x^1)^2} \crazy{\omega}{2}^2,
        \end{array}
       \right.
\renewcommand{\arraystretch}{1}
\end{equation}
\begin{equation}\label{eq:80}
          a_1 = - \displaystyle
\frac{1}{(x^1)^2} + \displaystyle \frac{1}{2 x^2 y^1}, \;\;\; a_2
= 0,
\end{equation}
\begin{equation}\label{eq:81}
\renewcommand{\arraystretch}{1.3}
\left\{
\begin{array}{ll}
 p_{11} = \displaystyle \frac{2}{(x^1)^3} +
\displaystyle \frac{x^1 y^2}{2 (x^2)^3 (y^1)^2},
\;\; p_{12} = -\displaystyle
\frac{1}{4 (x^2)^3 (y^1)^2}, \\
 q_{11} = -  \displaystyle
\frac{1}{2 x^2 (y^1)^2}, \;\; p_{2i} = q_{2i} =q_{12} = 0,
        \end{array}
        \right.
\renewcommand{\arraystretch}{1}
\end{equation}

\vspace*{-2mm}

\begin{equation}\label{eq:82}
p = -\displaystyle \frac{1}{8 (x^2)^3 (y^1)^2},
\;\; q = 0,
\end{equation}
\begin{equation}\label{eq:83}
\renewcommand{\arraystretch}{1.3}
\left\{
\begin{array}{ll}
b^1_{jkl} = b^2_{222}=  b^2_{122}= b^2_{212}= 0,   \;\;
b^2_{221} =  \displaystyle \frac{1}{2 (x^2)^3 (y^1)^2},\\
b^2_{111} = \displaystyle \frac{2 (y^2)^2(y^2 - x^1)}{(x^1)^4}
            + \displaystyle \frac{y^2 (x^1 + y^2)}{x^1 x^2 y^1}
            + \displaystyle \frac{x^1 y^2 (x^1y^2 + (x^2)^2)}{(x^2)^3
            (y^1)^2},\\
b^2_{112} = \displaystyle \frac{(2x^1-1) (4(x^2 y^1
            - (x^1)^2)}{4 (x^1)^4 x^2 y^1}, \\
b^2_{121} = - \displaystyle \frac{2(x^1)^2 (x^1 y^2
            - x^2 y^1) + (x^2)^2 y^1}{4 (x^1)^2 (x^2)^3 (y^1)^2}.
        \end{array}
        \right.
\renewcommand{\arraystretch}{1}
\end{equation}

It follows from  (80)--(83)
that equations (9), (13), (20), and (25) hold. Thus
the web (77) belongs to the classes ${\bf A_{21}}$
and ${\bf A_{22}}$.

Since by (83), $b^1_{222} \neq 0$,   the web (77) is not
transversally geodesic (and consequently is neither
hexagonal, nor Bol,  nor group web).  So this web
belongs to the class {\bf C}.
 By (81) and (82), the web (77) belongs to the classes
  ${\bf E_{21}}$ and ${\bf G_2}$.
}

\textbf{4. Examples of  nonextendable nonisoclinic webs}
$\boldsymbol{W} \boldsymbol{(}\boldsymbol{3}, \boldsymbol{2},
\boldsymbol{2}\boldsymbol{)}$
with $\boldsymbol{p} \boldsymbol{=} \boldsymbol{q} \boldsymbol{ \neq}
\boldsymbol{ 0}$ \textbf{(Class} ${\bf G_3).}$

\examp{\label{examp:5} A web $W(3,2,2)$ is given by the
equations
\begin{equation}\label{eq:84}
 u_3^1 = x^2 e^{x^1 y^1},\;\;\;
 u_3^2 = x^2 + y^2
\end{equation}
in a domain of ${\bf R}^4$ where $x^1, x^2, y^1 \neq 0$
(see [G 88], Example {\bf 8.1.27}, p. 391 or [G 92], p. 341
or [AS 92], Ch. 3, Problem {\bf 5}, p. 133)).

In this case using (48)--(51) and (32)--(37), we find that
\begin{equation}\label{eq:85}
\Gamma^2_{jk} = \Gamma^1_{i2}  = 0, \;\;
\Gamma^1_{11} = -\displaystyle
\frac{(1+x^1 y^1)e^{-x^1 y^1}}{x^1 x^2 y^1}, \;\;
\Gamma^1_{21} = \displaystyle \frac{1}{x^1 x^2 y^1},
\end{equation}
\begin{equation}\label{eq:86}
\omega^1_1 = - \displaystyle \frac{(1 + x^1 y^1) e^{-x^1 y^1}}{x^1 x^2 y^1}
\Bigl(\crazy{\omega}{1}^1 + \crazy{\omega}{2}^1\Bigr)
+ \displaystyle \frac{1}{x^1 x^2 y^1} \crazy{\omega}{1}^2, \;\;
 \omega_2^1 = \frac{1}{x^1 x^2 y^1} \crazy{\omega}{2}^1, \;\;
\end{equation}
\begin{equation}\label{eq:87}
  a_1 = 0, \;\; a_2 = - \displaystyle \frac{1}{x^1 x^2 y^1},
\end{equation}
\begin{equation}\label{eq:88}
p_{1i} = q_{ij} = 0, \;\; p_{22} = -\displaystyle
\frac{1 + x^1 y^1}{(x^1 x^2 y^1)^2}, \;\;
p_{21} = q_{21} = - \displaystyle \frac{e^{-x^1 y^1}}{(x^1 x^2 y^1)^2},
\end{equation}
 \begin{equation}\label{eq:89}
  p = q = \displaystyle \frac{e^{-x^1 y^1}}{2(x^1 x^2 y^1)^2},
\end{equation}
\begin{equation}\label{eq:90}
\renewcommand{\arraystretch}{1.3}
\left\{
\begin{array}{ll}
b^2_{jkl} = b^1_{jjj} =  b^1_{112}= b^1_{12i} = 0, \\
b^1_{211} =  -\displaystyle \frac{e^{-x^1 y^1}}{(x^1 x^2 y^1)^2},\;\;
2b^1_{212} = b^1_{221}= \displaystyle \frac{1 - x^1 y^1}{(x^1 x^2 y^1)^2}.
        \end{array}
        \right.
\renewcommand{\arraystretch}{1}
\end{equation}

It follows from  (87)--(90)
that equations (10), (14), (21), and (26) hold. Thus
the web (84) belongs to the classes ${\bf A_{31}}$
and ${\bf A_{32}}$.

Suppose that the web (84) is   transversally geodesic.
Then (38) implies that
$
b_{kl} = \frac{3}{4} b^i_{(ikl)}.
$
In particular, by (90), we must have
$
b_{12} =  \frac{1}{4} b^1_{211}.
$
Substituting this value of $b_{12}$ into (38), we find that
$
b^1_{(112)} = \frac{1}{6}  b^1_{211}.
$
But it follows from (90) that
$
b^1_{(112)} =  b^1_{211}.
$
Since $b^1_{211} \neq 0$ (see (90)), the web (84)
is not transversally geodesic (and consequently is neither
hexagonal, nor Bol,  nor group web).  So this web
belongs to the class {\bf C}.
 By (88) and (89), the web (84) belongs to the classes
 ${\bf E_{31}}$ and  ${\bf G_3}$.
}

\examp{\label{examp:6} Suppose that a web $W(3,2,2)$ is
given by the equations
\begin{equation}\label{eq:91}
 u_3^1 = x^1 + y^1 ,\;\;\;\;
 u_3^2 = - x^1 y^1 + x^2 y^2
\end{equation}
in a domain $x^2 \neq 0, y^2 \neq 0$
(see [G 88], Example {8.1.28}, p. 392,  [AS 92], Ch. {\bf 3},
Problem {\bf 4}, p. 133, [G 92], p. 341). In this case using
(48)--(51) and (32)--(37), we find that
\begin{equation}\label{eq:92}
\Gamma^1_{jk} = 0, \;\;  \Gamma^2_{11} =
 1-\displaystyle \frac{x^1 y^1}{x^2 y^2}, \;\;
\Gamma^2_{12} = -\displaystyle \frac{y^1}{x^2 y^2}, \;\;
\Gamma^2_{21} = - \displaystyle \frac{x^1}{x^2 y^2}, \;\;
\Gamma^2_{22} = -\displaystyle \frac{1}{x^2 y^2},
\end{equation}
\begin{equation}\label{eq:93}
 \renewcommand{\arraystretch}{1.3}
\left\{
\begin{array}{ll}
\omega^2_1 = \displaystyle  \Bigl(1-\displaystyle \frac{x^1 y^1}{x^2 y^2}\Bigr)
\Bigl(\crazy{\omega}{1}^1 + \crazy{\omega}{2}^1\Bigr)
-\displaystyle \frac{x^1}{x^2 y^2} \crazy{\omega}{1}^2
-\displaystyle \frac{y^1}{x^2 y^2} \crazy{\omega}{2}^2,
\\
\omega^2_2 = - \displaystyle \frac{1}{x^2 y^2}
\Bigl(\crazy{\omega}{1}^2 + \crazy{\omega}{2}^2\Bigr)
-\displaystyle \frac{y^1}{x^2 y^2} \crazy{\omega}{1}^1
-\displaystyle \frac{x^1}{x^2 y^2} \crazy{\omega}{2}^1, \;\;
 \omega^1_i =  0,
  \end{array}
\right.
\renewcommand{\arraystretch}{1}
\end{equation}
\begin{equation}\label{eq:94}
a_1 = \displaystyle \frac{x^1 - y^1}{x^2 y^2},\;\;\;\;
a_2 =0,
\end{equation}
\begin{equation}\label{eq:95}
 \renewcommand{\arraystretch}{1.3}
\left\{
\begin{array}{ll}
p_{2i} = q_{2i} = 0, &
 p_{12} = q_{12} = \displaystyle
\frac{y^1 - x^1}{(x^2 y^2)^2}, \\
 p_{11} = \displaystyle
\frac{x^2 y^1 - x^1 y^2 + (y^1)^2}{(x^2 y^2)^2}, &
 q_{11} = \displaystyle
\frac{- x^2 y^1 + x^1 y^2 - (x^1)^2}{(x^2 y^2)^2},
  \end{array}
\right.
\renewcommand{\arraystretch}{1}
\end{equation}
\begin{equation}\label{eq:96}
p = q = \displaystyle \frac{y^1 - x^1}{2(x^2 y^2)^2}
\end{equation}
\begin{equation}\label{eq:97}
\renewcommand{\arraystretch}{1.3}
\left\{
\begin{array}{ll}
b^1_{jkl} = b^2_{222}= b^2_{21i}=  0, \\
b^2_{111} =  \displaystyle
\frac{(x^1 - y^1)(x^2 y^2 - x^1 y^1)}{(x^2 y^2)^2},\;\;
b^2_{112} =  \displaystyle \frac{(y^1)^2 - x^1 y^1
                   + x^2 y^2}{(x^2 y^2)^2},\\
b^2_{121} =  \displaystyle \frac{-(x^1)^2 + x^1 y^1
                   - x^2 y^2}{(x^2 y^2)^2},\;\;
b^2_{122}= 2 b^2_{221} = \displaystyle \frac{y^1 - x^1}{(x^2
y^2)^2}.
        \end{array}
        \right.
\renewcommand{\arraystretch}{1}
\end{equation}

It follows from  (94)--(97)
that equations (9), (13), (20), and (25) hold. Thus
the web (91) belongs to the classes ${\bf A_{21}}$
and ${\bf A_{22}}$.

Since by (97), $b^2_{111} \neq 0$,   the web (91) is not
transversally geodesic.
% (and consequently is neither hexagonal, nor Bol,  nor group web).
 So this web
belongs to the class {\bf C}.
 By (95) and (96), the web (91) belongs to the classes
 ${\bf E_{23}}$ and ${\bf G_3}$.
 }

\examp{\label{examp:7} Suppose that a web $W(3,2,2)$ is
given by the equations
\begin{equation}\label{eq:98}
 u_3^1 = x^1 y^1  + x^2 (y^2)^2,\;\;\;\;
 u_3^2 =   x^2 + y^2
\end{equation}
in a domain $x^1 \neq 0, \; y^1 \neq 0$.
 In this case using (48)--(51) and (32)--(37), we find that
\begin{equation}\label{eq:99}
 \renewcommand{\arraystretch}{1.3}
\left\{
\begin{array}{ll}\Gamma^2_{jk} = 0, \;\;  \Gamma^1_{11} = -\displaystyle \frac{1}{x^1 y^1}, \;\;
\Gamma^1_{12} = \displaystyle \frac{2 x^2 y^2}{x^1 y^1}, \\
\Gamma^1_{21} = \displaystyle \frac{(y^2)^2}{x^1 y^1}, \;\;
\Gamma^1_{22} = -\displaystyle \frac{2 y^2(x^2 (y^2)^2 + x^1 y^1)}{x^1 y^1},
  \end{array}
\right.
\renewcommand{\arraystretch}{1}
\end{equation}
\begin{equation}\label{eq:100}
 \renewcommand{\arraystretch}{1.3}
\left\{
\begin{array}{ll}
\omega^1_1 =  - \displaystyle \frac{1}{x^1 y^1}
\Bigl(\crazy{\omega}{1}^1 + \crazy{\omega}{2}^1\Bigr)
+\displaystyle \frac{(y^2)^2}{x^1 y^1} \crazy{\omega}{1}^2
+\displaystyle \frac{2 x^2 y^2}{x^1 y^1} \crazy{\omega}{2}^2, \;\;
\omega^2_i =  0, \\
\omega^1_2 = -\displaystyle \frac{2 x^2 y^2}{x^1 y^1}
\Bigl(\crazy{\omega}{1}^1 + \crazy{\omega}{2}^2\Bigr)
+ \displaystyle \frac{2}{x^1 y^1}  \crazy{\omega}{1}^1
+ \displaystyle \frac{(y^2)^2}{x^1 y^1} \crazy{\omega}{2}^1,
   \end{array}
\right.
\renewcommand{\arraystretch}{1}
\end{equation}
\begin{equation}\label{eq:101}
 a_1 = 0, \;\; a_2 =
\displaystyle \frac{y^2 (-2 x^2 + y^2)}{x^1 y^1},
\end{equation}

\vspace*{-5mm}

\begin{equation}\label{eq:102}
\renewcommand{\arraystretch}{1.3}
\left\{
\begin{array}{ll}
 p_{1i} =
q_{1i} = 0, & p_{22} = \displaystyle \frac{(y^2)^2
[(y^2)^2 - 2 x^1 y^2 - 2 x^1 y^1]}{(x^1 y^1)^2},\\
 p_{21} =  q_{21} = \displaystyle \frac{(2x^2 -
y^2)y^2}{(x^1 y^1)^2},& q_{22} =
\displaystyle \frac{2 x^2 (y^2)^2 (y^2 - 2 x^2)}{(x^1 y^1)^2},
         \end{array}
\right.
\renewcommand{\arraystretch}{1}
\end{equation}
\begin{equation}\label{eq:103}
 p = q = - \displaystyle \frac{1}{p_{21}},
\end{equation}
\begin{equation}\label{eq:104}
\renewcommand{\arraystretch}{1.3}
\left\{
\begin{array}{ll}
b^2_{jkl} = b^1_{1kl}= b^1_{211} = 0, &
b^1_{222} =  \displaystyle \frac{2 u_3^1 [(y^2)^2 (2x^2 - y^2) + x^1 y^1]}{(x^1 y^1)^2}\\
b^1_{212} =  \displaystyle \frac{2 x^2 (y^2)^2 (y^2 - 2x^2)}{(x^1
y^1)^2},&
b^1_{221}= \displaystyle \frac{(y^2)^3 (y^2 - 2 x^2) - 2 x^1 y^1 y^2}{(x^1 y^1)^2}.
        \end{array}
        \right.
\renewcommand{\arraystretch}{1}
\end{equation}

It follows from  (101)--(104)
that equations (10), (14), (21), and (26) hold. Thus
the web (98) belongs to the classes ${\bf A_{31}}$
and ${\bf A_{32}}$.

The web (98) is not  transversally geodesic since
$b^1_{222} \neq 0$ (see (104)). Thus the web (98) belongs to the class {\bf C}.
 By (102) and (103), the web (98) belongs to the classes
  ${\bf E_{32}}$ and  ${\bf G_3}$.
}

\examp{\label{examp:8} {\it Polynomial $3$-webs.} We will
call a web $W(3,2,2)$ {\it polynomial} if it is defined by
the equations
\begin{equation}\label{eq:105}
u_3^i = x^i + y^i + c_{jk}^i x^j y^k,
\end{equation}
where $c_{jk}^i = \mbox{const.}$ and
\begin{eqnarray*}
\renewcommand{\arraystretch}{1.3}
 \Delta_1 \!\!\!\!& = \!\!\!\! & (1+c_{1i}^1 y^i ) (1+ c_{2j}^2 y^j ) - c_{1i}^2
 c_{2j}^1 y^i y^j \neq 0, \\
 \Delta_2 \!\!\!\!& = \!\!\!\!& (1+c_{i1}^1 x^i ) (1+ c_{j2}^2 x^j ) - c_{i1}^2 c_{j2}^1
 x^i x^j \neq 0
\renewcommand{\arraystretch}{1}
\end{eqnarray*}
(see [G 88], Example {\bf 8.4.4}, p. 430).

For a polynomial web, using (48), (50) and (6), we
obtain
\begin{equation}\label{eq:106}
\renewcommand{\arraystretch}{1.3}
\left\{
\begin{array}{ll}
         a_1 =  \displaystyle
\frac{1}{\Delta_1 \Delta_2} \Bigl[&\!\!\!\!\!c_{21}^2 - c_{12}^2
+ (c_{11}^2 c_{22}^2 -
                     c_{12}^2 c_{21}^2)
                     (y^1 - x^1) \\
&\!\!\!\!\!+ (c_{12}^1 c_{11}^2 -
                     c_{11}^1 c_{12}^2) x^1
 +  (c_{22}^1 c_{11}^2 - c_{21}^1 c_{12}^2)x^2 \\
&\!\!\!\!\!+ (c_{11}^1   c_{21}^2 - c_{21}^1 c_{11}^2) y^1 +
(c_{12}^1 c_{21}^2 - c_{22}^1 c_{11}^2) y^2\Bigr], \\
         a_2 =   \displaystyle
\frac{1}{\Delta_1 \Delta_2} \Bigl[&\!\!\!\!\!c_{12}^1 - c_{21}^1
+ (c_{11}^1
                     c_{22}^1 - c_{12}^1 c_{21}^1)
                     (y^2 - x^2) \\
&\!\!\!\!\!+ (c_{22}^1 c_{11}^2 -  c_{21}^1 c_{12}^2)x^1
            +  (c_{22}^1 c_{21}^2 -   c_{21}^1 c_{22}^2 ) x^2 \\
&\!\!\!\!\!+ (c_{12}^1  c_{21}^2 - c_{22}^1 c_{11}^2 ) y^1
      +  (c_{12}^1 c_{22}^2 - c_{22}^1 c_{12}^2) y^2\Bigr].
        \end{array}
        \right.
\renewcommand{\arraystretch}{1.3}
\end{equation}
Differentiating equations (106) and applying (4), we can calculate
the values of $p_{ij}$ and $q_{ij}$. Their expressions
obtained by a lengthy calculations are rather long.
However,  it appeared that
$$
p_{12} - p_{21} = q_{12} - q_{21} = \displaystyle \frac{1}{\Delta_1^2 \Delta_2^2}
\Bigl[(c^1_{11} + c^2_{21})(c^1_{12} - c^1_{21})
- (c^2_{22} + c^1_{12})(c^2_{21} - c^2_{12})\Bigr],
$$
i.e., we have
\begin{equation}\label{eq:107}
p = q = \displaystyle \frac{1}{2\Delta_1^2 \Delta_2^2}
\Bigl[(c^1_{11} + c^2_{21})(c^1_{12} - c^1_{21})
- (c^2_{22} + c^1_{12})(c^2_{21} - c^2_{12})\Bigr].
\end{equation}
Equation (107) shows that the same condition
\begin{equation}\label{eq:108}
(c_{22}^2 + c_{12}^1 )(c_{12}^2 - c_{21}^2) + (c_{11}^1 +
c_{21}^2 ) (c_{12}^1 - c_{21}^1) \neq 0
\end{equation}
is  sufficient for both, $p \neq 0$ and $q \neq 0$. Therefore, a
polynomial web (105) satisfying condition (108) is nonisoclinic.

As a result,  the polynomial
3-webs (105) cannot be extended to a nonisoclinic
$W(4,2,2)$, and hence it belongs to the class
${\bf G_3}$.

Note that  condition (108) will not be satisfied
 for example if one of the following conditions holds:
\begin{itemize}
\item $c^1_{12} = c^1_{21}, \;\; c^2_{12} = c^2_{21}$.
\item $c^2_{22} = - c^1_{12}, \;\; c^2_{21} = - c^1_{11}$.
\item $c^2_{22} =- c^1_{21}, \;\; c^2_{12} = - c^2_{11}$.
\item $c^2_{22} = - c^1_{12} = -c^2_{21}$.
\item $c^1_{11} = - c^2_{12} = - c^2_{21}$.
\end{itemize}
In these cases  further investigation is necessary to determine
whether a web in question is isoclinic or nonisoclinic, and if it
is nonisoclinic, to which of the classes ${\bf G_1-G_4}$ it belongs.
}

\examp{\label{examp:9} Take a particular case of the web
(105) by taking $c^1_{jk} = 0$:
\begin{equation}\label{eq:109}
u_3^1 = x^1 + y^1, \;\; u_3^2 = x^2 + y^2 + c_{jk}^2 x^j y^k,\;\;
 c_{jk}^2 = \mbox{const}.
\end{equation}
where
$$
\Delta_1 = 1 +  c^2_{21} y^1 + c^2_{22} y^2 \neq 0, \;\;
\Delta_2 = 2 +  c^2_{12} x^1 + c^2_{22} x^2 \neq 0.
$$
In this case using (48)--(51) and (32)--(37), we find that
\begin{equation}\label{eq:110}
 \renewcommand{\arraystretch}{1.3}
\left\{
\begin{array}{ll}
\Gamma^1_{jk} = 0, \;\;
\Gamma^2_{22}
    = \displaystyle \frac{c^2_{22}}{\Delta_1 \Delta_2}, \;\;
\Gamma^2_{12} = \displaystyle \frac{\alpha y^1 - c^2_{12}}{\Delta_1 \Delta_2}, \;\;
\Gamma^2_{21} = \displaystyle \frac{\alpha x^1 - c^2_{21}}{\Delta_1
\Delta_2},\\
\Gamma^2_{11}
= - \displaystyle \frac{\Delta_1 (\alpha x^2 + c^2_{11}) + (c^2_{11} y^1 +  c^2_{12}y^2)
   (\alpha x^1 - c^2_{21})}{\Delta_1 \Delta_2}, \\
     \end{array}
\right.
\renewcommand{\arraystretch}{1}
\end{equation}
\begin{equation}\label{eq:111}
 a_1 =  -\displaystyle
\frac{\beta}{\Delta_1 \Delta_2}, \;\; a_2 = 0,
\end{equation}
\begin{equation}\label{eq:112}
\renewcommand{\arraystretch}{1.3}
\left\{
\begin{array}{ll}
p_{2i} = q_{2i} = 0, & p_{11} = \displaystyle
\frac{-\alpha \Delta_1^2 \Delta_2 + \beta(c^2_{12} - \alpha y^1)}{\Delta_1^3
\Delta_2^2},\\
 p_{12} =  q_{12} =  \displaystyle \frac{c^2_{22}\beta}{\Delta_1^2 \Delta_2^2 },
&
q_{11} = \displaystyle
\frac{\alpha \Delta_1 \Delta_2^2 + \beta (c^2_{21} - \alpha x^1)}{\Delta_1^2 \Delta_2^3},
        \end{array}
        \right.
\renewcommand{\arraystretch}{1.3}
\end{equation}
\begin{equation}\label{eq:113}
 p = q =  \displaystyle \frac{c^2_{22}\beta}{2\Delta_1^2 \Delta_2^2 },
\end{equation}
where
$$
\alpha = c^2_{11} c^2_{22} -  c^2_{12} c^2_{21}, \;\;
\beta = \alpha (x^1 - y^1) + c^2_{12} - c^2_{21}.
        $$

In this case the condition (108) becomes $$ c^2_{22} \neq 0. $$
Thus for this web $p \neq 0$ and $q \neq 0$
if $c^2_{22} \neq 0$. Conditions (112) and (113)
prove that the web (109) belongs to the classes
  ${\bf E_{23}}$ and ${\bf G_3}$.
Note that the condition $c^2_{22} = 0$ is necessary
and sufficient for the web (109) to be isoclinic.

By (111), (112) and (9) and (20), the  web (109) belongs
to the class  ${\bf A_{21}}$. It is interesting to
find out for which values of $c^2_{jk}$ the web (109) belongs
to the classes ${\bf C}$ or ${\bf D}$ and to the class ${\bf A_{22}}$.
   }

{\examp{\label{examp:10} Take another particular case of the web
(105) by taking $c^1_{11} = c^1_{21} = c^1_{22} = c^2_{11} =
c^2_{12} = c^2_{22} = 0$ and $c^1_{12} = c^2_{21} = 1$:
\begin{equation}\label{eq:114}
u_3^1 = x^1 + y^1 + x^1 y^2,  \;\; u_3^2 = x^2 + y^2 + x^2 y^1
\end{equation}
in a domain
$$
\Delta_1 = (1 + y^1) (1 + y^2) \neq 0, \;\;
\Delta_2 = 1 - x^1 x^2 \neq 0.
$$
 i.e., $y^i \neq - 1, \;\; x^1 x^2 \neq 1$.
In this case using (48)--(51) and (32)--(37), we find that
\begin{equation}\label{eq:115}
\renewcommand{\arraystretch}{1.3}
\left\{
\begin{array}{ll}
\Gamma^1_{11} = - x^2 \Gamma^1_{12}
=  \displaystyle \frac{x^2 y^2}{\Delta_1\Delta_2},&
\Gamma^2_{1i} = \Gamma^1_{2i}  = 0,\\
\Gamma^2_{22} = - x^1 \Gamma^2_{21}
= \displaystyle \frac{x^1}{(1+ y^1)\Delta_2}, &
      \end{array}
       \right.
\renewcommand{\arraystretch}{1}
\end{equation}
\begin{equation}\label{eq:116}
\renewcommand{\arraystretch}{1.3}
\left\{
\begin{array}{ll}
\omega^1_1 = a_2 x^2 (\crazy{\omega}{1}^1 + \crazy{\omega}{2}^1)
- a_2 \crazy{\omega}{2}^2, &
 \omega_1^2 = - a_1 \crazy{\omega}{1}^2, \\
  \omega^2_2 =   a_1 x^1 (\crazy{\omega}{1}^2 +
 \crazy{\omega}{2}^2) - a_1  \crazy{\omega}{2}^1, &
 \omega_2^1 = - a_2 \crazy{\omega}{1}^1,
      \end{array}
       \right.
\renewcommand{\arraystretch}{1}
\end{equation}
\begin{equation}\label{eq:117}
  a_1 =  \displaystyle
\frac{1}{(1 + y^1) \Delta_2}, \;\;
 a_2 = \displaystyle
\frac{y^2}{\Delta_1 \Delta_2},
\end{equation}
\begin{equation}\label{eq:118}
\renewcommand{\arraystretch}{1.3}
\left\{
\begin{array}{ll}
p_{12} = q_{12} =a_1 (a_1 x^1 + a_2), & q_{11}
= - a_1 (a_2 x^2 + a_1), \\
p_{21} = q_{21} = a_2 (a_2 x^2 + a_1), &
q_{22} = - a_{2} (a_1 x^1 + a_2),  \\
p_{11} = p_{22} = 0. &
        \end{array}
        \right.
\renewcommand{\arraystretch}{1.3}
\end{equation}
\begin{equation}\label{eq:119}
 p = q = \displaystyle
\frac{a_1^2 x^1 - a_2^2 x^2}{2},
\end{equation}
\begin{equation}\label{eq:120}
\renewcommand{\arraystretch}{1.3}
\left\{
\begin{array}{ll}
b^i_{jjj} = b^1_{12i}= b^1_{2i1}= b^2_{i12} = b^2_{2i1}  =0, & \\
b^1_{112} = - \displaystyle \frac{x^2}{2 (1 + y^2)^2 \Delta_2^2},
\;\;\;\;\;\;\;\;\;\;\;\;\;\;\;\;\;
b^1_{212} = - \displaystyle \frac{x^1 + y^1 + x^1 y^2 + 1}{(1 + y^2) \Delta_1 \Delta_2^2}, \\
b^2_{121} = - \displaystyle \frac{x^2 + y^2 + x^2 y^1 + 1}{(1 + y^1) \Delta_1
\Delta_2^2}, \;\;\;\;\;\;\;\;
b^2_{122} =  \displaystyle \frac{x^1 + y^1 + x^1 y^2 + 1}{(1 + y^1) \Delta_1
\Delta_2^2}.
        \end{array}
        \right.
\renewcommand{\arraystretch}{1}
\end{equation}

By (118) and (119), the web (114) belongs to the classes
  ${\bf E_{131}}$ and ${\bf G_3}$.
Using (117) and (118), it is easy to check that
condition (8) does not hold. Thus
the web (114) belongs to the class {\bf B}.

Suppose that the web (114) is  transversally geodesic. Then
it follows from (38) and (120) that
$b_{22} = \frac{3}{4} b^i_{(i22)} = \frac{1}{4} b^1_{212} $.
As a result, by (38) and (120), we have
$b^1_{(122)} = \frac{1}{3} b_{22} = \frac{1}{12} b^1_{212}$.
On the other hand, $b^1_{(122)} = \frac{1}{3} b^1_{212}$.
Since by (120), $b^1_{212} \neq 0$,
we came to a contradiction.
This contradiction proves that the web (114)
is not transversally geodesic, i.e., this web
belongs to the class {\bf C}.
 }

\examp{\label{examp:11}
 Consider the web  defined by the equations
\begin{equation}\label{eq:121}
u_3^1 =  x^1 +  y^1 +  \displaystyle \frac{1}{2} (x^1)^2 y^2,
\;\;\; u_3^2 =  x^2 + y^2 +  \displaystyle \frac{1}{2} x^1 (y^2)^2
\end{equation}
 in a domain of ${\bf R}^4$ where $\Delta = 1 + x^1 y^2  \neq 0$,
 i.e., $x^1 y^2  \neq -1$.

In this case using (48)--(51) and (32)--(37), we find that
\begin{equation}\label{eq:122}
 \Gamma^1_{12} = -\displaystyle \frac{x^1}{\Delta^2}, \;\;\;
 \Gamma^2_{12} = - \displaystyle \frac{y^2}{\Delta^2}, \;\;
 \Gamma^1_{i1} = \Gamma^1_{22} = \Gamma^2_{i1} = \Gamma^2_{22} = 0,
\end{equation}
\begin{equation}\label{eq:123}
\omega^1_1 =  -\displaystyle \frac{x^1}{\Delta^2}
\crazy{\omega}{2}^2, \;\; \omega_2^1 = -\displaystyle \frac{x^1}{\Delta^2}
\crazy{\omega}{1}^1, \;\;  \omega_1^2
= -\displaystyle \frac{y^2}{\Delta^2}
\crazy{\omega}{2}^2, \;\; \omega_2^2
= -\displaystyle \frac{y^2}{\Delta^2}
\crazy{\omega}{1}^1,
\end{equation}
\begin{equation}\label{eq:124}
a_1 =- \displaystyle \frac{y^2}{\Delta^2}, \;\; a_2 =
\displaystyle \frac{x^1}{\Delta^2},
\end{equation}
\begin{equation}\label{eq:125}
\renewcommand{\arraystretch}{1.3}
\left\{
\begin{array}{ll}
p_{i2} = 0, & q_{i1} = 0, \\
p_{11} = \displaystyle \frac{2 (y^2)^2}{\Delta^4}, &
q_{22} = -\displaystyle \frac{2 (x^1)^2}{\Delta^4}, \\
p_{21} = - q_{12} = \displaystyle \frac{1 - x^1 y^2}{\Delta^4},
&
\end{array}
\right.
\renewcommand{\arraystretch}{1}
\end{equation}
\begin{equation}\label{eq:126}
 p = q = \displaystyle \frac{x^1 y^2 - 1}{2 \Delta^4},
\end{equation}
\begin{equation}\label{eq:127}
\renewcommand{\arraystretch}{1.3}
\left\{
\begin{array}{ll}
b^i_{jjj} = b^1_{12i}= b^1_{2i1}=  b^2_{2i1} = b^2_{12i} = 0, & \\
b^1_{112} = b^2_{112} = - b^2_{212}= \displaystyle \frac{x^1 y^2 -
1}{\Delta^4},\;\;
b^1_{212} = - \displaystyle \frac{2 (x^1)^2}{\Delta^3}.
        \end{array}
        \right.
\renewcommand{\arraystretch}{1}
\end{equation}

By (125) and (126), the web (121) belongs to the classes
${\bf G_3}$ and ${\bf E_{111}}$. Recall that for webs of the class
${\bf E_{111}}$, both components $a_1$ and
$a_2$ are covariantly  constant on the 2-dimensional distribution defined
by the equations
  $\crazy{\omega}{1}^1 = \crazy{\omega}{2}^2  = 0$, and they satisfy
      the exterior quadratic equation $ \nabla a_1 \wedge \crazy{\omega}{1}^1
      = \nabla a_2 \wedge \crazy{\omega}{2}^2$ which is equivalent
      to the condition $p_{21} = - q_{12}$.
We establish now a geometric meaning of the latter condition for
the web (121).

Consider the 3-dimensional distribution defined by the equation
 $ a_1  \crazy{\omega}{1}^1  - a_2   \crazy{\omega}{2}^2 = 0$.
 Taking exterior derivative of this equation and applying
 (123) and (124), we find that
 $$
      \nabla a_1 \wedge \crazy{\omega}{1}^1
      - \nabla a_2 \wedge \crazy{\omega}{2}^2 = 0.
 $$
It follows that {\em for the web $(121)$,  the condition
 $p_{21} = - q_{12}$ is necessary and sufficient for
 the distribution  $ a_1  \crazy{\omega}{1}^1  - a_2
  \crazy{\omega}{2}^2 = 0$ to be integrable.}

Note also that  by (123) and (124), we have
$$
d(\crazy{\omega}{1}^1 + \crazy{\omega}{2}^1) = 0, \;\;
d(\crazy{\omega}{1}^2 + \crazy{\omega}{2}^2) = 0.
$$
This means that {\em for the web $(121)$, the $3$-dimensional distributions defined
by the equations $\crazy{\omega}{1}^1 + \crazy{\omega}{2}^1 = 0$
and $\crazy{\omega}{1}^2 + \crazy{\omega}{2}^2 = 0$
are integrable}.

Using (124) and (125), it is easy to check that
the left-hand side of condition (8) is
$\displaystyle \frac{2 x^1 (x^2)^2 + x^1 y^2 (1 - x^1 y^2)}{\Delta^8} \neq 0$,
i.e., condition (8) does not hold. Thus
the web (121) belongs to the class \textbf{B}.

Suppose that the web (121) is transversally geodesic.
Then it follows from (38) and (127) that
$b_{12} = \frac{1}{4} b^i_{(i22)} = 0$.
As a result, by (38) and (127), we have
$b^1_{(112)} = \frac{2}{3} b_{12} = 0$.
On the other hand $b^1_{(112)}
= \frac{1}{3} b^1_{112}$.
Since by (127), $b^1_{112} \neq 0$,
we come to a contradiction.
This contradiction proves that the web (121)
is not transversally geodesic, i.e., this web
belongs to the class \textbf{C}.
}

\examp{\label{examp:12}
 Consider the web  defined by the equations
\begin{equation}\label{eq:128}
u_3^1 =  (x^1 +  y^1)(x^2 - y^2), \;\;\; u_3^2 =  (x^1 - y^1)
(x^2 + y^2)
\end{equation}
 in a domain of ${\bf R}^4$ where $\Delta  = x^1 y^2  + x^2 y^1 \neq 0$.

In this case using (48)--(51) and (32)--(37), we find that
\begin{equation}\label{eq:129}
\renewcommand{\arraystretch}{1.3}
\left\{
\begin{array}{ll}
\Gamma^1_{11} = - \Gamma^2_{11} = \displaystyle \frac{u_3^2}{2\Delta^2},
\;\;\;\;
\Gamma^1_{22} = -\Gamma^2_{22} =  -\displaystyle \frac{u_3^1}{2\Delta^2},
\\
\Gamma^1_{12} = - \Gamma^1_{21} =  -\Gamma^2_{12} = \Gamma^2_{21}
=\displaystyle \frac{\rho}{2\Delta^2}, \;\;\;
       \end{array}
       \right.
\renewcommand{\arraystretch}{1}
\end{equation}
\begin{equation}\label{eq:130}
\renewcommand{\arraystretch}{1.3}
\left\{
\begin{array}{ll}
\omega^1_1 = - \omega_1^2 = \displaystyle \frac{1}{2\Delta^2}
\bigl[u_3^2(\crazy{\omega}{1}^1 + \crazy{\omega}{2}^1)
- \rho (\crazy{\omega}{1}^2 - \crazy{\omega}{2}^2)\bigr], \\
\omega^2_2 = - \omega_2^1 = \displaystyle \frac{1}{2\Delta^2}
\bigl[\rho (\crazy{\omega}{2}^1 - \crazy{\omega}{1}^1)
+ u^1_3(\crazy{\omega}{1}^1 + \crazy{\omega}{2}^1)\bigr],
       \end{array}
       \right.
\renewcommand{\arraystretch}{1}
\end{equation}
\begin{equation}\label{eq:131}
a_1 = a_2 =- \displaystyle \frac{\rho}{2\Delta^2},
\end{equation}
\begin{equation}\label{eq:132}
\renewcommand{\arraystretch}{1.3}
\left\{
\begin{array}{ll}
 p_{11} = p_{21} = \displaystyle
\frac{\beta (x^1 - y^1) - \alpha (x^2 +  y^2)}{2\Delta^4},  \\
 p_{12} = p_{22} = \displaystyle \frac{- \beta (x^1 +
y^1) + \alpha (x^2 -  y^2)}{2\Delta^4}, \\
 q_{11} = q_{21} = \displaystyle
\frac{\gamma (x^1 - y^1) + \delta  (x^2 +  y^2)}{2\Delta^4},\\
 q_{12} = q_{22} = \displaystyle
\frac{\gamma (x^1 +y^1) + \delta(x^2 - y^2)}{2 \Delta^4},
\end{array}
\right.
\renewcommand{\arraystretch}{1}
\end{equation}
\begin{equation}\label{eq:133}
 p = q =
\displaystyle \frac{\rho \sigma}{\Delta^4},
\end{equation}
\begin{equation}\label{eq:134}
\renewcommand{\arraystretch}{1.3}
\left\{
\begin{array}{ll}
b^2_{i11} = - b^1_{i11} = \displaystyle \frac{u_3^2 \rho}{\Delta^4},
\;\; b^1_{222} = - b^2_{222} = - b^i_{122} = \displaystyle \frac{u_3^1
\rho}{\Delta^4},\\
b^1_{i12} = - b^1_{i21} = - b^2_{i12} = b^2_{i21}
=  \displaystyle \frac{\sigma^2
- 2 [(x^1 x^2)^2 + (y^1 y^2)^2]}{\Delta^4},
\end{array}
\right.
\renewcommand{\arraystretch}{1}
\end{equation}
where
$$
\renewcommand{\arraystretch}{1.3}
\left\{
\begin{array}{ll}
\alpha =(x^1)^2 y^2 - x^1 x^2 y^1 - 2
(y^1)^2 y^2, &
\beta = (x^2)^2 y^1 - x^1 x^2 y^2 - 2 y^1 (y^2)^2, \\
\gamma = x^1 (y^2)^2 -  x^2 y^1 y^2 - 2 x^1 (x^2)^2, &
\delta = x^2  (y^1)^2 - x^1 y^1 y^2 - 2 (x^1)^2 x^2,\\
\rho = x^1 x^2 + y^1 y^2, & \sigma = x^1 y^2 - x^2 y^1.
\end{array}
\right.
\renewcommand{\arraystretch}{1}
$$

By  (133),  the web (128) belongs to the class
${\bf G_3}$.
Equations (131) and (132) show that
conditions (11) and (22) hold. Thus,
the web (128) belongs to the class ${\bf A_{131}}$.
It is easy to check that the curvature tensor
of this web defined by equations (134) satisfies both
equations (27). As a result,
the web (128) belongs to the class ${\bf A_{132}}$.

It follows from  (134) that $b^1_{222} \neq  0$,
and by Corollary 8a),  this  proves that the web (121)
is not transversally geodesic, i.e., this web
belongs to the class {\bf C}.
}

\examp{\label{examp:13}
 Consider the web  defined by the equations
\begin{equation}\label{eq:135}
u_3^1 = x^2 y^2 e^{x^1 y^1}, \;\;\; u_3^2 =  x^2 + y^2
\end{equation}
in a domain of ${\bf R}^4$ where $x^i \neq 0, \; y^i \neq 0$.
In this case using (48)--(51) and (32)--(37), we find that
\begin{equation}\label{eq:136}
 \renewcommand{\arraystretch}{1.3}
\left\{
\begin{array}{ll}
\Gamma^1_{12} = \displaystyle \frac{1}{x^1 y^1 y^2}, \;\;
\Gamma^1_{21} = \displaystyle \frac{1}{ x^1 x^2 y^1},\;\;
\Gamma^1_{22} = - \displaystyle \frac{e^{x^1 y^1}}{x^1 y^1}, \\
\Gamma^1_{11}
    = -\displaystyle \frac{(1 + x^1 y^1) e^{-x^1 y^1}}{x^1 x^2  y^1 y^2}, \;\;
\Gamma^2_{jk} = 0,
     \end{array}
\right.
 \renewcommand{\arraystretch}{1}
\end{equation}
\begin{equation}\label{eq:137}
 \renewcommand{\arraystretch}{1.3}
\left\{
\begin{array}{ll}
\omega_1^1 =  \displaystyle \frac{1}{x^1 x^2 y^1 y^2}
\Bigl[- (1 + x^1 y^1) e^{-x^1 y^1}
 (\crazy{\omega}{1}^1 + \crazy{\omega}{2}^1)
+  y^2 \crazy{\omega}{1}^2 + x^2\crazy{\omega}{2}^2\Bigr], \\
\omega_2^1 = \displaystyle \frac{1}{x^1 x^2 y^1 y^2}
\Bigl[x^2\crazy{\omega}{1}^1 +  y^2\crazy{\omega}{2}^1
- x^2 y^2 e^{x^1 y^1} (\crazy{\omega}{1}^2
+  \crazy{\omega}{2}^2)\Bigr], \;\; \omega_i^2 = 0,
     \end{array}
\right.
 \renewcommand{\arraystretch}{1}
\end{equation}
\begin{equation}\label{eq:138}
a_1 = 0, \;\;\;\; a_2 =  \displaystyle \frac{x^2 - y^2}{x^1 x^2 y^1
y^2},
\end{equation}
\begin{equation}\label{eq:139}
\renewcommand{\arraystretch}{1.3}
\left\{
\begin{array}{ll}
p_{1i} = q_{1i} = 0, & p_{21} = q_{21} =
 \displaystyle \frac{(y^2 - x^2)e^{-x^1 y^1}}{(x^1 x^2 y^1
y^2)^2}, \\
 p_{22} =  \displaystyle \frac{x^2 - y^2 + x^1 y^1 y^2}{(x^1 x^2 y^1)^2 y^2}, &
 q_{22} =  \displaystyle \frac{x^2 - y^2 - x^1 x^2 y^1}{(x^1 y^1 y^2)^2 x^2},
   \end{array}
\right.
\renewcommand{\arraystretch}{1}
\end{equation}
\begin{equation}\label{eq:140}
p = q = \displaystyle \frac{(x^2 - y^2) e^{- x^1
y^1}}{2(x^1 x^2 y^1 y^2)^2}.
\end{equation}
\begin{equation}\label{eq:141}
\renewcommand{\arraystretch}{1.3}
\left\{
\begin{array}{ll}
b^2_{jkl} = b^1_{11i} =   b^1_{12i} =  0, & \\
b^1_{211} = \displaystyle \frac{(x^2 - y^2) (2 - x^1 y^1) e^{-x^1 y^1}}{(x^1 x^2 y^1 y^2)^2}, &
 b^1_{212} = \displaystyle \frac{2(y^2 - x^2) + x^1 x^2 y^1}{(x^1 y^1 y^2)^2 x^2}, \\
b^1_{222} = \displaystyle \frac{(x^2 - y^2) (2 - x^1 y^1) e^{x^1 y^1}}{(x^1 y^1)^2 x^2 y^2},
&
b^1_{221} = \displaystyle \frac{2(y^2 - x^2) + x^1 y^1 y^2}{(x^1 x^2 y^1)^2
y^2}.
   \end{array}
\right.
\renewcommand{\arraystretch}{1}
\end{equation}

By (139) and (140),  the web (135) belongs to the classes
 ${\bf E_{32}}$ and ${\bf G_3}$.
Equa\-tions (138), (139), and (141) show that
conditions (10), (21), and (26) hold. Thus,
the web (135) belongs to the classes ${\bf A_{31}}$
and ${\bf A_{32}}$.

It follows from (141) that $b^1_{222} \neq 0$, and
as a result, according to Corollary 8a),
 the web (135) is not transversally geodesic, i.e., this web
belongs to the class {\bf C}.
 }

\examp{\label{examp:14}
 Consider the web  defined by the equations
\begin{equation}\label{eq:142}
u_3^1 = y^2 e^{x^1 y^1}, \;\;\; u_3^2 =  x^2 + y^2
\end{equation}
in a domain of ${\bf R}^4$ where $x^1 \neq 0, \; y^i \neq 0$
(see [AS 92], Ch. 3, Problem {\bf 5}, p. 133).
In this case using (48)--(51) and (32)--(37), we find that
\begin{equation}\label{eq:143}
 \renewcommand{\arraystretch}{1.3}
\left\{
\begin{array}{ll}
\Gamma^1_{11}
  = -\displaystyle \frac{(1 + x^1 y^1) e^{-x^1 y^1}}{x^1  y^1  y^2},\;\;
    \Gamma^1_{12} = \displaystyle \frac{1}{x^1 y^1 y^2}, \\
\Gamma^1_{2i} = 0,  \;\;
\Gamma^2_{jk} = 0,
     \end{array}
\right.
 \renewcommand{\arraystretch}{1}
\end{equation}
\begin{equation}\label{eq:144}
 \renewcommand{\arraystretch}{1.3}
\left\{
\begin{array}{ll}
\omega_1^1 =  \displaystyle \frac{1}{x^1 y^1 y^2}
\Bigl[- (1 + x^1 y^1) e^{-x^1 y^1}
 (\crazy{\omega}{1}^1 + \crazy{\omega}{2}^1)
+  \crazy{\omega}{2}^2\Bigr], \\
\omega_2^1 = \displaystyle \frac{1}{x^1  y^1 y^2}
\crazy{\omega}{1}^1, \;\; \omega_i^2 = 0,
     \end{array}
\right.
 \renewcommand{\arraystretch}{1}
\end{equation}
\begin{equation}\label{eq:145}
a_1 = 0, \;\;\;\; a_2 =  -\displaystyle \frac{1}{x^1 y^1 y^2},
\end{equation}
\begin{equation}\label{eq:146}
\renewcommand{\arraystretch}{1.3}
\left\{
\begin{array}{ll}
p_{1i} = q_{1i} = 0, & p_{21} = q_{21} =
 \displaystyle \frac{e^{-x^1 y^1}}{(x^1 y^1 y^2)^2}, \\
 p_{22} =  0, &
 q_{22} =  -\displaystyle \frac{1}{(x^1 y^1 y^2)^2},
   \end{array}
\right.
\renewcommand{\arraystretch}{1}
\end{equation}
\begin{equation}\label{eq:147}
p = q = -\displaystyle \frac{e^{- x^1 y^1}}{2(x^1 y^1 y^2)^2}.
\end{equation}
\begin{equation}\label{eq:148}
\renewcommand{\arraystretch}{1.3}
\left\{
\begin{array}{ll}
b^2_{jkl} = b^1_{11i} =   b^1_{12i} =  b^1_{22i}= 0, & \\
b^1_{211} = \displaystyle \frac{e^{-x^1 y^1} (3 + 2 x^1 y^1)}{(x^1 y^1 y^2)^2}, \;\;\;
 b^1_{212} = \displaystyle \frac{x^1 y^1 - 1}{(x^1 y^1 y^2)^2}.
   \end{array}
\right.
\renewcommand{\arraystretch}{1}
\end{equation}

By (146) and (147),  the web (142) belongs to the classes
 ${\bf E_{321}}$ and ${\bf G_3}$.
Equations (145), (146), and (148) show that
conditions (10), (21), and (26) hold. Thus,
the web (142) belongs to the classes ${\bf A_{31}}$
and ${\bf A_{32}}$.

It follows from (148) that $b^2_{jkl} = 0$, and
as a result, according to Corollary 8b),
 the web (142) is not transversally geodesic, i.e., this web
belongs to the class {\bf C}.

Note that in  [AS 92], Ch. 3,
Problem {\bf 5}, p. 133, it is wrongly indicated that the web
(142) is transversally geodesic.
 }

\examp{\label{examp:15}
 Consider the web  defined by the equations
\begin{equation}\label{eq:149}
u_3^1 = e^{x^1 y^1} + x^2 y^2 , \;\;\; u_3^2 =  x^2 + y^2
\end{equation}
in a domain of ${\bf R}^4$ where $x^1, \, y^2 \neq 0$
(see [AS 92], Ch. 3, Problem {\bf 8}, p. 133).
In this case using (48)--(51) and (32)--(37), we find that
\begin{equation}\label{eq:150}
 \renewcommand{\arraystretch}{1.3}
\left\{
\begin{array}{ll}
\Gamma^1_{11} = - A  e^{-x^1 y^1},\;\;
\Gamma^1_{12} = A x^2 e^{-x^1 y^1}, \;\;
\Gamma^1_{21} = A y^2 e^{-x^1 y^1},\\
\Gamma^1_{22} = - A x^2 y^2 e^{-x^1 y^1} - 1,\;\;
\Gamma^2_{jk} = 0,
     \end{array}
\right.
 \renewcommand{\arraystretch}{1}
\end{equation}
\begin{equation}\label{eq:151}
 \renewcommand{\arraystretch}{1.3}
\left\{
\begin{array}{ll}
\omega_1^1 =   A  e^{-x^1 y^1}[- (\crazy{\omega}{1}^1 + \crazy{\omega}{2}^1)
+   y^2  \crazy{\omega}{1}^2 +   x^2 \crazy{\omega}{2}^2], \;\; \omega_i^2 = 0,\\
\omega_2^1 =  A e^{-x^1 y^1}[ (x^2 y^2  + A^{-1} e^{x^1 y^1}) [- (\crazy{\omega}{1}^2
+  \crazy{\omega}{2}^2) + x^2 \crazy{\omega}{1}^1
+    y^2  \crazy{\omega}{2}^1],
      \end{array}
\right.
 \renewcommand{\arraystretch}{1}
\end{equation}
\begin{equation}\label{eq:152}
a_1 = 0, \;\; a_2 = A(x^2 - y^2)e^{-x^1 y^1},
\end{equation}
\begin{equation}\label{eq:153}
\renewcommand{\arraystretch}{1.3}
\left\{
\begin{array}{ll}
p_{1i} = q_{1i} = 0, \;\;
 p_{22} = - y^2 p_{21} +   A e^{-x^1 y^1}, \;\;
 q_{22} = - x^2 p_{21}, \\
 p_{21} = q_{21} =
 (y^2 - x^2) Be^{-2x^1 y^1},
\end{array}
\right.
\renewcommand{\arraystretch}{1}
\end{equation}
\begin{equation}\label{eq:154}
p = q = \displaystyle \frac{(x^2 -
y^2)B  e^{- 2 x^1 y^1}}{2},
\end{equation}
\begin{equation}\label{eq:155}
\renewcommand{\arraystretch}{1.3}
\left\{
\begin{array}{ll}
b^2_{jkl} = b^1_{11i} =   b^1_{12i} = 0,  \\
b^1_{211} = (x^2 - y^2) (B + A^2) e^{-2x^1 y^1} \\
 b^1_{212} = A e^{-x^1 y^1} -  B (x^2)^2 e^{-2x^1 y^1}, \\
b^1_{221} = - A e^{-x^1 y^1} - (x^2 - y^2) y^2 (B + A^2) e^{-2x^1 y^1},\\
b^1_{222} = \displaystyle \frac{1}{2}
\bigl(4 x^2 - 3y^2\bigr)\bigl(A  e^{-x^1 y^1} + x^2 y^2 B  e^{-2x^1 y^1}\bigr),
\end{array}
\right.
\renewcommand{\arraystretch}{1}
\end{equation}
where
$$
A = 1 + \displaystyle \frac{1}{x^1 y^1}, \;\;
B = A + \displaystyle \frac{1}{(x^1 y^1)^2}.
$$

By (153) and (154),  the web (149) belongs to the classes
 ${\bf E_{32}}$ and ${\bf G_3}$. \linebreak
Equations (152), (153), and (155) show that
conditions (10), (21), and (26) hold. Thus,
the web (142) belongs to the classes ${\bf A_{31}}$
and ${\bf A_{32}}$.

It follows from (155) that $b^1_{222} \neq 0$, and
as a result according to Corollary 8a),
 the web (142) is not transversally geodesic, i.e., this web
belongs to the class {\bf C}.
}

 \textbf{5. Examples of  nonextendable nonisoclinic
webs} $\boldsymbol{W} \boldsymbol{(}\boldsymbol{3},
\boldsymbol{2}, \boldsymbol{2}\boldsymbol{)}$ with
$\boldsymbol{p} \boldsymbol{\neq} \boldsymbol{0}, \;
\boldsymbol{q} \boldsymbol{\neq} \boldsymbol{0},\;
\boldsymbol{p} \boldsymbol{\neq} \boldsymbol{q},$
\textbf{and
condition (43) not held (Class} ${\bf G_4).}$

\examp{\label{examp:16} A web $W(3,2,r)$ is given by
\begin{equation}\label{eq:156}
 u_3^1 = \displaystyle
\frac{1}{6} (x^1 + y^1)^3 + \displaystyle \frac{1}{2}
\Bigl[(x^1)^2 + (y^1)^2 + 2x^2 y^2\Bigr] ,\;\;\;
 u_3^2 = x^2 + y^2
\end{equation}
in a domain where
$$
 \Delta_1 =\displaystyle  \frac{1}{2} \Bigl(x^1 + y^1\Bigr)^2 +
x^1 \neq 0,\;\;\;\;
 \Delta_2 =\displaystyle \frac{1}{2} \Bigl(x^1
+ y^1\Bigr)^2 + y^1 \neq 0
$$
(see [G 88], Ch. 8, Example {\bf 8.1.29}, p. 392
and [G 92]).

In this case using (48)--(51) and (32)--(37), we find that
\begin{equation}\label{eq:157}
 \renewcommand{\arraystretch}{1.3}
\left\{
\begin{array}{ll}
\Gamma^1_{11} = - \displaystyle\frac{\alpha}{\Delta_1 \Delta_2}, &
\Gamma^1_{12} =  \displaystyle\frac{\alpha x^2}{\Delta_1 \Delta_2},\\
\Gamma^1_{21} = \displaystyle\frac{\alpha y^2}{\Delta_1 \Delta_2},&
\Gamma^1_{22} =  - \Bigl[\displaystyle\frac{\alpha x^2 y^2}{\Delta_1 \Delta_2}
+ 1\Bigr],\;\; \Gamma^2_{jk} = 0,
     \end{array}
\right.
 \renewcommand{\arraystretch}{1}
\end{equation}
\begin{equation}\label{eq:158}
 \renewcommand{\arraystretch}{1.3}
\left\{
\begin{array}{ll}
\omega_1^1 =  \displaystyle\frac{\alpha}{\Delta_1 \Delta_2}
 \Bigl[- (\crazy{\omega}{1}^1 + \crazy{\omega}{2}^1)
+   y^2  \crazy{\omega}{1}^2 +   x^2 \crazy{\omega}{2}^2\Bigr], \;\; \omega_i^2 = 0,\\
\omega_2^1 = \displaystyle\frac{\alpha}{\Delta_1 \Delta_2}\Bigl[-\Bigl(x^2 y^2  +
\displaystyle\frac{\Delta_1 \Delta_2}{\alpha}\Bigr) (\crazy{\omega}{1}^2
+  \crazy{\omega}{2}^2) + x^2 \crazy{\omega}{1}^1
+    y^2  \crazy{\omega}{2}^1\Bigr],
      \end{array}
\right.
 \renewcommand{\arraystretch}{1}
\end{equation}
\begin{equation}\label{eq:159}
 a_1 = 0 ,\;\; a_2 = -
\displaystyle \frac{\alpha \beta}{\Delta_1 \Delta_2},
\end{equation}
\begin{equation}\label{eq:160}
\renewcommand{\arraystretch}{1.3}
\left\{
\begin{array}{lll}
p_{1i} = 0, &p_{21} = \displaystyle
 \frac{\beta A}{\Delta_1^3
\Delta_2^2}, & p_{22} = - y^2 p_{21}
- \displaystyle\frac{\alpha}{\Delta_1
\Delta_2},\\ q_{1i} = 0, &q_{21} = \displaystyle \frac{\beta
B}{\Delta_1^2 \Delta_2^3}, &  q_{22} = - x^2 q_{21}
+ \displaystyle\frac{\alpha}{\Delta_1 \Delta_2},
        \end{array}
        \right.
\renewcommand{\arraystretch}{1}
\end{equation}
\begin{equation}\label{eq:161}
         p  = - \displaystyle
\frac{\beta A}{2\Delta_1^3 \Delta_2^2}, \;\;
         q   = - \displaystyle
\frac{\beta B}{2\Delta_1^2 \Delta_2^3},
\end{equation}
\begin{equation}\label{eq:162}
\renewcommand{\arraystretch}{1.3}
\left\{
\begin{array}{ll}
b^1_{111} =  \displaystyle \frac{x^1 - y^1}{\Delta_1^3 \Delta_2^3}
\Bigl(\displaystyle -\frac{3}{4} \alpha^4 - \frac{3}{2} \alpha^3 - \alpha^2
+ x^1 y^1\Bigr),  \\
b^1_{112} = \displaystyle \frac{(x^1 - y^1) x^2}{\Delta_1^3 \Delta_2^3}
\Bigl(\displaystyle \frac{3}{4} \alpha^4 - \frac{1}{2} \alpha^3 - \alpha^2
- x^1 y^1\Bigr),  \\
 b^1_{121} = b^1_{211} = \displaystyle \frac{1}{\Delta_1^3 \Delta_2^3}
\Bigl[ (y^2 \Delta_2 - x^2 \Delta_1)\bigl(\Delta_1 \Delta_2 - \alpha^2
(\Delta_1 + \Delta_2)\bigr) \\
\hspace*{34mm} +  \alpha (x^2 \Delta_1^2 - y^2 \Delta_2^2)\Bigr], \\
b^1_{122} =  \displaystyle \frac{\alpha (\alpha + 1) x^2 y^2}{\Delta_1^3 \Delta_2^3}
\Bigl[\displaystyle \frac{1}{4} \alpha^4 + \frac{1}{2} \alpha^3
+ \alpha (\alpha + 1) (x^1- y^1) + x^1 y^1\Bigr],  \\
b^1_{212} =  \displaystyle \frac{1}{\Delta_1^3 \Delta_2^3}
\Bigl[ \alpha^2 (\alpha + 1) (\Delta_2 y^2 \beta - \Delta_1 (x^2)^2)
+ y^2 \beta \Delta_2^2 (2\alpha - \frac{1}{2}  \Delta_1)\\
\hspace*{23mm} + (x^2)^2 \Delta_1^2 (1 - \alpha)\Bigr], \\
b^1_{221} =  \displaystyle \frac{\alpha}{4\Delta_1^3 \Delta_2^3}
\Biggl\{ y^2 (x^2 \Delta_1 - y^2  \Delta_2)
\biggl[\frac{4 x^1 y^1}{(x^2 \Delta_1 - y^2  \Delta_2) \alpha}
- \alpha^2 (3\alpha + 2)\biggr] \\
\hspace*{24mm} + 4x^2 (y^2 \Delta_2^2 - x^2   \Delta_1^2) - 4\Delta_1^2 \Delta_2^2\Biggl\},\\
b^1_{222} =  \displaystyle \frac{\alpha}{\Delta_1^3 \Delta_2^3}
\Bigl[- \beta \Delta_1^2 \Delta_2^2 + x^2 y^2 (\alpha + 1)\bigl(\Delta_1 \Delta_2
- \alpha^2 (x^2 \Delta_1 + y^2 \Delta_2)\bigr)\\
\hspace*{23mm} -  x^2 y^2 (x^2 \Delta_1^2 + y^2 \Delta_2^2)\Bigr],
\;\; b^2_{jkl} =  0.
\end{array}
\right.
\renewcommand{\arraystretch}{1}
\end{equation}
where
$$
\renewcommand{\arraystretch}{1.3}
\left\{
\begin{array}{ll}
\alpha = x^1 + y^1, & \beta = x^2 - y^2, \\
A = \displaystyle \frac{3}{4} \alpha^4 + \alpha^3 + (y^1)^2, &
B = \displaystyle \frac{3}{4} \alpha^4 + \alpha^3 + (x^1)^2.
\end{array}
\right.
\renewcommand{\arraystretch}{1}
$$

By (160),  the web (156) belongs to the class
 ${\bf E_{3}}$.
Equations (159), (160), and (162) show that
conditions (10) and (21) hold but condition (26)
does not hold. Thus,
the web (156) belongs to the class ${\bf A_{31}}$
and does not belong to the class ${\bf A_{32}}$.

It follows from (162) that $b^1_{222} \neq 0$, and
as a result, according to Corollary 8a),
 the web (156) is not transversally geodesic, i.e., this web
belongs to the class {\bf C}.

We now prove that {\em the web $(156)$ cannot
be expanded to a web $W (4, 2, 2)$ of maximum $2$-rank.}
Thus we have to prove that
the relations (43) do not hold for the web (149).

It follows from (36) and (151) that
$$
d p = p \omega_1^1 + \crazy{p}{1}_i^{ } \crazy{\omega}{1}^i +
 \crazy{p}{2}_i^{ } \crazy{\omega}{2}^i ,\;\;\;
 d q = q \omega_1^1 + \crazy{q}{1}_i^{ } \crazy{\omega}{1}^i +
 \crazy{q}{2}_i^{ } \crazy{\omega}{2}^i,
$$
where  $\omega_1^1$ is defined by (158).
On the other hand, we can find $d p$ and $d q$ by
differentiating their expressions (161). Comparing the results, we find that
\begin{equation}\label{eq:163}
 \renewcommand{\arraystretch}{1.3}
\left\{
\begin{array}{ll}
\crazy{p}{1}_1 =  -\displaystyle\frac{3 \alpha^2 (\alpha + 1)
\beta}{2\Delta_1^4 \Delta_2^2} + \displaystyle\frac{3 A (\alpha +
1)  \beta}{2\Delta_1^5 \Delta_2^2} + \displaystyle\frac{A \alpha
\beta}{2\Delta_1^4 \Delta_2^3}, & \crazy{p}{1}_2 = - y^2
\crazy{p}{1}_1 -\displaystyle\frac{A}{2\Delta_1^3 \Delta_2^2}, \\
\crazy{p}{2}_1 =  -\displaystyle\frac{[3\alpha^2 (\alpha + 1) + 2
y^1] \beta}{2\Delta_1^3 \Delta_2^3} + \displaystyle\frac{ A
\alpha \beta}{\Delta_1^4 \Delta_2^3} + \displaystyle\frac{A
(\alpha + 1) \beta}{\Delta_1^3 \Delta_2^4}, & \crazy{p}{2}_2 = -
x^2 \crazy{p}{2}_1 + \displaystyle\frac{A}{2\Delta_1^3
\Delta_2^2}, \\
 \crazy{q}{1}_1 =  -\displaystyle\frac{[3 \alpha^2 (\alpha +
 1)+ 2 x^1]\beta }{2\Delta_1^3 \Delta_2^3}
+ \displaystyle\frac{B \alpha \beta}{\Delta_1^3 \Delta_2^4}
+ \displaystyle\frac{B (\alpha + 1) \beta }{\Delta_1^4
\Delta_2^3},&
\crazy{q}{1}_2 = - y^2 \crazy{q}{1}_1  -\displaystyle\frac{B}{2\Delta_1^2 \Delta_2^3}, \\
\crazy{q}{2}_1 =  -\displaystyle\frac{3 \alpha^2 (\alpha +
1) \beta}{2\Delta_1^2 \Delta_2^4}
+ \displaystyle\frac{3B (\alpha + 1) \beta}{2\Delta_1^2 \Delta_2^5}
+ \displaystyle\frac{B \alpha \beta}{2\Delta_1^3 \Delta_2^4}, &
\crazy{q}{2}_2 = - x^2 \crazy{q}{2}_1
+ \displaystyle\frac{B}{2\Delta_1^2 \Delta_2^3}.
      \end{array}
\right.
 \renewcommand{\arraystretch}{1}
\end{equation}

Let us assume that equation (43) holds for $i = 1$. Since
by (159), $a_1 = 0$, this means that
\begin{equation}\label{eq:164}
 \crazy{q}{ }_{ }^{ } (\crazy{q}{ }_{ }^{ }
 \crazy{p}{1}_1^{ } - \crazy{p}{ }_{ }^{ }
 \crazy{q}{1}_1^{ }) -\crazy{p}{ }_{ }^{ }
 (\crazy{q}{ }_{ }^{ }\crazy{p}{2}_1^{ } -
 \crazy{p}{ }_{ }^{ } \crazy{q}{2}_1^{ }) = 0.
\end{equation}
By (163) and (164), the left-hand side  LHS of
 equation (43) taken for $i = 2$ is
 $$
\mbox{{\rm LHS}} = \crazy{q}{ }_{ }^{ } (\crazy{q}{ }_{ }^{ }
 \crazy{p}{1}_2^{ } - \crazy{p}{ }_{ }^{ }
 \crazy{q}{1}_2^{ }) -\crazy{p}{ }_{ }^{ }
 (\crazy{q}{ }_{ }^{ }\crazy{p}{2}_2^{ } -
 \crazy{p}{ }_{ }^{ } \crazy{q}{2}_2^{ })
 = \beta p (q \crazy{p}{2}_1 - p \crazy{q}{2}_1).
$$
Applying (163),  from the last equation we find that
\begin{equation}\label{eq:165}
 \renewcommand{\arraystretch}{1.3}
\begin{array}{ll}
\mbox{{\rm LHS}} = \displaystyle \frac{\beta^3 p}{4\Delta_1^6 \Delta_2^7}
\Biggl\{AB \Bigl[(\alpha + 1) x^1 - \alpha y^1 + \frac{1}{2} \alpha^2\Bigr]
\\
\hspace*{27mm}
+ \Delta_1 \Delta_2 \Bigl[3 \alpha^3 (\alpha + 1) (x^1 - y^1)
+ 2 B y^1\Bigr]\Biggr\}.
      \end{array}
 \renewcommand{\arraystretch}{1}
\end{equation}

On the other hand, by (159) and (161), for $i = 2$
the right-hand side RHS of (43) is
 \begin{equation}\label{eq:166}
\mbox{{\rm RHS}}
= \displaystyle \frac{\beta^3 p}{4\Delta_1^6 \Delta_2^7}
\Biggl\{AB (-A \Delta_2 + B \Delta_1)\Biggr\}.
\end{equation}
It is easy to see that in the curl brackets of expression (165) the
highest degree of $x^1$ is 11 while in those of (166)
the highest degree of $x^1$ is 14. This
proves that the 2nd equation of (43) fails.
Thus, by Theorem 7, part e), the web (149) cannot
be expanded to a web $W (4, 2, 2)$ of maximum 2-rank.
As a result, this web belongs to the class
${\bf G_4}$.
}

\examp{\label{examp:17}
 Consider the web  defined by
\begin{equation}\label{eq:167}
u_3^1 =  x^1 +  y^1 +  \displaystyle \frac{1}{2} (x^1)^2 y^2,
\;\;\; u_3^2 =  x^2 + y^2 +  \displaystyle \frac{1}{2} (x^2)^2 y^1
\end{equation}
 in a domain of ${\bf R}^4$ where $\Delta_1 = (1 + x^1 y^2)(1 + x^2 y^1) \neq 0$
and $\Delta_2 = \displaystyle \frac{1}{4} \Bigl[4 - (x^1 x^2)^2 \Bigr] \neq 0$,
i.e., $ x^1 y^2 \neq -1, \; x^2 y^1 \neq -1, \; x^1 x^2 \neq \pm 2$.

In this case using (48)--(51) and (32)--(37), we find that
\begin{equation}\label{eq:168}
 \renewcommand{\arraystretch}{1.3}
\left\{
\begin{array}{ll}
\Gamma^1_{11} =  \displaystyle\frac{\beta x^1 (x^2)^2}{2\Delta_1 \Delta_2}, &
\Gamma^2_{22} =  \displaystyle\frac{\alpha (x^1)^2 x^2}{2\Delta_1\Delta_2},\\
\Gamma^1_{12} = - \displaystyle\frac{\beta x^1}{\Delta_1 \Delta_2},&
\Gamma^2_{21} = - \displaystyle\frac{\alpha x^2}{\Delta_1 \Delta_2},\\
 \Gamma^1_{2i} =  = 0, &\Gamma^2_{1i} =  0,
     \end{array}
\right.
 \renewcommand{\arraystretch}{1}
\end{equation}

\begin{equation}\label{eq:169}
 \renewcommand{\arraystretch}{1.3}
\left\{
\begin{array}{ll}
\omega_1^1 =  \displaystyle\frac{\beta x^1 (x^2)^2}{2\Delta_1 \Delta_2}
  (\crazy{\omega}{1}^1 + \crazy{\omega}{2}^1)
- \displaystyle\frac{\beta x^1}{\Delta_1 \Delta_2} \crazy{\omega}{2}^2,&
\omega_2^1 =  -  \displaystyle\frac{\beta x^1}{\Delta_1 \Delta_2}
\crazy{\omega}{1}^1, \\
\omega_2^2 =  \displaystyle\frac{\alpha (x^1)^2 x^2}{2\Delta_1 \Delta_2}
  (\crazy{\omega}{1}^2 + \crazy{\omega}{2}^2)
-  \displaystyle\frac{\alpha x^2}{\Delta_1 \Delta_2}\crazy{\omega}{2}^1,&
\omega_1^2 =  -  \displaystyle\frac{\alpha x^2}{\Delta_1 \Delta_2}
\crazy{\omega}{1}^2, \\
      \end{array}
\right.
 \renewcommand{\arraystretch}{1}
\end{equation}

\begin{equation}\label{eq:170}
a_1 = \displaystyle \frac{\alpha x^2}{\Delta_1 \Delta_2},
\;\; a_2 = \displaystyle \frac{\beta x^1}{\Delta_1
\Delta_2},
\end{equation}
\begin{equation}\label{eq:171}
\renewcommand{\arraystretch}{1.3}
\left\{
\begin{array}{ll}
p_{11} = 0, &
q_{11} = - \displaystyle
\frac{A(x^2)^2}{2\Delta_1 \Delta_2^2 \beta}, \\
p_{12} = \displaystyle
\frac{1}{\Delta_2 \beta^3} +
\frac{(x^1 x^2)^2}{2 \Delta_2^2 \beta^2}
+ \frac{x^1 x^2}{\Delta_1 \Delta_2}, &
q_{12} =  \displaystyle
\frac{Bx^1 x^2}{2 \Delta_1 \Delta_2^2 \beta}, \\
 p_{21} = \displaystyle
\frac{1}{\Delta_2 \alpha^3} +
\frac{(x^1 x^2)^2}{2 \Delta_2^2 \alpha^2}
+ \frac{x^1 x^2}{\Delta_1 \Delta_2}, &
q_{21} = \displaystyle
\frac{Ax^1 x^2}{2 \Delta_1 \Delta_2^2 \alpha}, \\
 p_{22} = 0, &
q_{22} = - \displaystyle
\frac{B(x^1)^2}{2 \Delta_1 \Delta_2^2 \alpha},
        \end{array}
        \right.
\renewcommand{\arraystretch}{1}
\end{equation}
\begin{equation}\label{eq:172}
 p = q + \displaystyle
\frac{1}{2\Delta_2}\Biggl(\displaystyle \frac{1}{\beta^3}
- \displaystyle \frac{1}{\alpha^3}\Biggr), \;\;
 q = \displaystyle
\frac{(x^1 x^2)^2}{4\Delta_2^2}\Biggl(\displaystyle \frac{1}{\beta^2}
- \displaystyle \frac{1}{\alpha^2}\Biggr),
\end{equation}
\begin{equation}\label{eq:173}
\renewcommand{\arraystretch}{1.3}
\left\{
\begin{array}{ll}
b^1_{111} = \displaystyle \frac{(x^2)^2}{2 \alpha^3
\Delta_2^2}\Bigl[\Delta_2 + \displaystyle \frac{\alpha (x^1 x^2)^2}{4}\Bigr],
& b^1_{121} = -\displaystyle \frac{x^1 x^2}{4 \Delta_1
 \Delta_2^2}, \\
b^1_{112} = \displaystyle \frac{1}{4 \alpha^3
\Delta_2^2}\Bigl[ -2\Delta_2
+  \alpha (x^1 x^2)^2\Bigr], &
 b^1_{211} = \displaystyle \frac{x^1 x^2}{8 \alpha \Delta_1
 \Delta_2^2} \Bigl[4\alpha - \beta x^1 (x^2)^3\Bigr], \\
 b^1_{212} = -\displaystyle \frac{(x^1)^2}{4\alpha \Delta_2^2}
 \Bigl[\alpha x^1 x^2 (1 + x^2) + 2 \beta (2 + (x^2)^2)\Bigr], &
b^1_{i22} =   b^1_{221}  = 0, \\
 b^2_{222} = \displaystyle \frac{(x^1)^2}{8 \beta^3 \Delta_2^2}
\Bigl[4\Delta_2 + \beta (x^1 x^2)^2\Bigr],
  &
 b^2_{121} = -\displaystyle \frac{A (x^2)^2}{2\beta \Delta_1 \Delta_2^2},
 \\
 b^2_{122} = \displaystyle \frac{x^1 x^2}{2\beta \Delta_1 \Delta_2}
 \Bigl[2 \alpha x^1 x^2 + \beta (1 - x^2)\Bigr], &
  b^2_{212} = \displaystyle \frac{x^1 (x^2)^2}{2 \Delta_1 \Delta_2^2},
 \\
 b^2_{221} = -   \displaystyle \frac{(x^1)^2}{4 \beta^3 \Delta_2^2}
\Bigl[4\Delta_2 + \beta (x^1 x^2)^2\Bigr],
 &
 b^2_{i11} =  b^2_{112} =  0,
        \end{array}
        \right.
\renewcommand{\arraystretch}{1}
\end{equation}
where
$$
\renewcommand{\arraystretch}{1.3}
\left\{
\begin{array}{ll}
 \alpha = 1 + x^1 y^2, & \beta = 1 + x^2 y^1,\\
   A = 2 \alpha + \beta x^1 x^2,
 &B = 2 \beta + \alpha x^1 x^2.
       \end{array}
        \right.
\renewcommand{\arraystretch}{1}
$$

It follows from equations (170) and (171) that conditions (8) do
not hold. In fact, by (171), the first and the third
terms of the left-hand side of the first equation of (8) vanish. Up to a common
factor, the middle term of the first equation of (8) is
$$
2 \alpha \beta \Delta_2 (\alpha^3 + \beta^3) + (x^1)^2  (x^2)^2
\alpha^2 \beta^2 (\alpha^2 + \beta^2) + 2 x^1 x^2 \Bigl(1 - \displaystyle
\frac{(x^1)^2 (x^2)^2}{4}\Bigr).
$$
The degree of the 2nd term of this expression (with respect
to $x^i$ and $y^j$) is 16, and it is higher than the degrees
14 and 6 of two other terms. This term is
$$
(x^1)^2  (x^2)^2 (1 + x^1 y^2)^2 (1 + x^2 y^1)^2 (2 + (x^1)^2 (y^2)^2 +
(x^2 y^1)^2 + 2 x^1 y^2 + 2 x^2 y^1).
$$
The highest degree terms of the last expression are
$$
(x^1)^4  (x^2)^4  (y^1)^2 (y^2)^2 ((x^1)^2 (y^2)^2 +
(x^2 y^1)^2)
$$
do not vanish and have no similar terms with two other terms of
the left-hand side of the first equation of (8).
 Thus the web (167) belongs to the class {\bf B}.

Equations (171) prove that  the web (167) belongs to the class ${\bf E_{13}}$.

It follows from (36) that
$$
d p = p (\omega_1^1  + \omega_2^2)+ \crazy{p}{1}_i^{ } \crazy{\omega}{1}^i +
 \crazy{p}{2}_i^{ } \crazy{\omega}{2}^i ,\;\;\;
 d q = q (\omega_1^1 + \omega_2^2) + \crazy{q}{1}_i^{ } \crazy{\omega}{1}^i +
 \crazy{q}{2}_i^{ } \crazy{\omega}{2}^i,
$$
where  $\omega_1^1$ and $\omega_2^2$ are defined by (169).
On the other hand, we can find $d p$ and $d q$ by
differentiating their expressions (172). Comparing the results, we find that
\begin{equation}\label{eq:174}
 \renewcommand{\arraystretch}{1.3}
\left\{
\begin{array}{ll}
\crazy{q}{1}_1  = \displaystyle\frac{x^1 (x^2)^2 (4 + (x^1 x^2)^2)}{4\alpha\Delta_2^3}
\Bigl(\displaystyle \frac{1}{\beta^2} - \displaystyle
\frac{1}{\alpha^2}\Bigr) + \displaystyle \frac{(x^1 x^2)^2 y^2}{\alpha^4
\Delta_2^2},\\
\crazy{q}{1}_2  = \displaystyle\frac{(x^1)^2 x^2 (4 + (x^1 x^2)^2)}{4\beta\Delta_2^3}
\Bigl(\displaystyle \frac{1}{\beta^2} - \displaystyle
\frac{1}{\alpha^2}\Bigr) - \displaystyle \frac{(x^1 x^2)^2 y^1}{\beta^4
\Delta_2^2},\\
\crazy{q}{2}_1  = \displaystyle\frac{(x^1)^2 (x^2)^3}{4\alpha \Delta_2^2}
\Bigl[\displaystyle \frac{3 x^1 y^1}{\alpha^2 \Delta_2} - \displaystyle
\frac{1}{\beta}\Bigl(\displaystyle \frac{x^1 x^2}{\beta}
+ \displaystyle \frac{2}{\alpha}\Bigl)\Bigr], \;\;
\crazy{p}{1}_1 = \crazy{q}{1}_1  + \displaystyle\frac{3 y^2}{2\alpha^5\Delta_2},\\
\crazy{q}{2}_2  = \displaystyle\frac{(x^1)^3 (x^2)^2}{4\beta \Delta_1^2\Delta_2^3}
\Bigl[x^1 y^1 \Bigl(\alpha^2 + \beta^2\Bigr) + 2
\Delta_1\Bigr], \;\;
\crazy{p}{1}_2 = \crazy{q}{1}_2  - \displaystyle\frac{3 y^1}{2\beta^5\Delta_2},\\
\crazy{p}{2}_1 = \crazy{q}{2}_1  - 2K - (x^2)^2 L,\;\;
 \crazy{p}{2}_2 = \crazy{q}{2}_2  + (x^1)^2 K + 2 L,
       \end{array}
\right.
 \renewcommand{\arraystretch}{1}
\end{equation}
where
$$
K = \displaystyle \frac{x^2}{4\beta\Delta_2^2}
\Biggl(\displaystyle \frac{1}{\alpha^3} + \displaystyle
\frac{2}{\beta^3}\Biggr), \;\;
L =  \displaystyle \frac{x^1}{4\alpha\Delta_2^2}
\Biggl(\displaystyle \frac{2}{\alpha^3} + \displaystyle
\frac{1}{\beta^3}\Biggr).
$$

Consider the first of two equations (43). Substitute
the values of $p, q, a_1$,
and $\crazy{p}{\alpha}_1, \;  \crazy{q}{\alpha}_1$
from (172), (170), and (174) into this equation, and
collect all similar terms. It turns out that
the highest (53rd) degree term is $- 48 (x^1 x^2)^{18} (y^1)^9
(y^2)^8 \neq 0$. Thus, the web (167) belongs to the class
${\bf G_4}$.

Finally, we prove that the web (167) is not transversally
geodesic. Suppose that it is transversally
geodesic.   By (173), the component $b^1_{112}$ of the curvature
tensor of the web (167) vanishes,  $b^1_{112} = 0$.
Since this web  is transversally
geodesic, by (38), we have   $b^1_{112} = \frac{1}{3} b_{11}$.
But as we noted earlier, $b_{11} = \frac{3}{4} b^k_{(k11)}$.
By (173), we find that $b_{11} = \frac{1}{4} (3b^1_{111}
+ b^2_{121})$. A straightforward calculation shows that
$b_{11}$ is proportional to a polynomial in $x^i$ and
$y^j$ whose  highest degree term is $S (x^1)^3 (x^2)^4 y^2 (y^1)^2$
with $S \neq 0$. Thus $b_{11} \neq 0$, and
consequently $b^1_{112} \neq 0$. This contradiction proves that
the web (167) belongs to the class {\bf C}.
}

\examp{\label{examp:18}
 Consider the web  defined by  the equations
\begin{equation}\label{eq:175}
u_3^1 =  x^1 +  y^1 +  x^1 y^2, \;\;\; u_3^2 =  x^1 y^1 + x^2 y^2
\end{equation}
in a domain of ${\bf R}^4$ where $\Delta_1 = y^2 (1+y^2) \neq 0$,
and $\Delta_2 = x^2 - (x^1)^2 \neq 0$, i.e.,
 $y^2  \neq 0, -1$ and $(x^1)^2 \neq x^2$.

In this case using (48)--(51) and (32)--(37), we find that

\begin{equation}\label{eq:176}
 \renewcommand{\arraystretch}{1.3}
\left\{
\begin{array}{lll}
\Gamma^1_{11} =  \displaystyle\frac{x^1 y^2}{\Delta_1 \Delta_2}, &
\Gamma^1_{12} =  -\displaystyle\frac{x^2 y^2}{\Delta_1 \Delta_2}, &
 \Gamma^2_{12} =  \displaystyle\frac{x^1 y^2 +
y^1}{\Delta_1\Delta_2},\\
\Gamma^2_{22} = -  \displaystyle\frac{1}{y^2 \Delta_2}, &
\Gamma^2_{11} =  -\displaystyle\frac{x^1 y^1 + x^2
y^2}{\Delta_1\Delta_2}, &
 \Gamma^2_{21} =  \displaystyle\frac{x^1}{y^2\Delta_2},\;\; \Gamma^1_{2i} = 0,
     \end{array}
\right.
 \renewcommand{\arraystretch}{1}
\end{equation}
\begin{equation}\label{eq:177}
 \renewcommand{\arraystretch}{1.3}
\left\{
\begin{array}{lll}
\omega_1^1 =  \displaystyle\frac{y^2}{\Delta_1 \Delta_2}
  \Bigl[x^1 (\crazy{\omega}{1}^1 + \crazy{\omega}{2}^1)
- x^2  \crazy{\omega}{2}^2\Bigr], \;\;\;\;
\omega_2^1 =  -\displaystyle\frac{x^2 y^2}{\Delta_1 \Delta_2}
\crazy{\omega}{1}^1,\\
\omega_1^2 = -\displaystyle\frac{1}{\Delta_1 \Delta_2}
  \Bigl[(x^1 y^1 + x^2 y^2) (\crazy{\omega}{1}^1 + \crazy{\omega}{2}^1)
+ x^1 (1 + y^2)  \crazy{\omega}{1}^2 + (x^1 y^2 + y^1)   \crazy{\omega}{2}^2\Bigr], \\
\omega_2^2 = \displaystyle\frac{1}{\Delta_1 \Delta_2}\Bigr[
   (x^1 y^2 + y^1) \crazy{\omega}{1}^1 - (1 + y^2)
   (\crazy{\omega}{1}^2 + \crazy{\omega}{2}^2) + x^1 (1 + y^2)  \crazy{\omega}{2}^1\Bigr],
     \end{array}
\right.
 \renewcommand{\arraystretch}{1}
\end{equation}
\begin{equation}\label{eq:178}
a_1 = \displaystyle \frac{y^1 - x^1}{\Delta_1 \Delta_2}, \;\;
a_2 = \displaystyle \frac{x^2 y^2}{\Delta_1 \Delta_2 },
\end{equation}
\begin{equation}\label{eq:179}
 \renewcommand{\arraystretch}{1.3}
\left\{
\begin{array}{ll}
p_{11} = \displaystyle
\frac{1}{\Delta_1^2 \Delta_2^2} \Bigl[(x^1 y^1 + x^2 y^2)(1 - x^2 y^2)
- x^1 y^2 (2 x^1 - 3 y^1)  + (y^1)^2\Bigr], \\
q_{11} = \displaystyle
\frac{1}{\Delta_1^2 \Delta_2^2} \Bigl[(x^1 y^1 + x^2 y^2)(1 + x^2  +y^2)
- (x^1)^2 (1 +  y^2)\Bigr], \\
      p_{12} = \displaystyle
\frac{1}{\Delta_1^2 \Delta_2 y^2}
\Bigl[x^1 - y^1 - x^1 x^2 y^2\Bigr], \;\;\;
 \\
  p_{21} = \displaystyle
\frac{1}{\Delta_1^2\Delta_2^2}
\Bigl[x^1 y^2 (x^2 y^2 +  x^1 y^1 + x^1 - y^1) - x^2 y^1 y^2\Bigr], \\
 q_{12} =  \displaystyle
\frac{1}{\Delta_1^2 \Delta_2^2}
\Bigl[(x^1 + y^1)y^2 + (1 + y^2)(x^1 + y^1 - x^1 y^2
- x^1 x^2 y^2)\Bigr], \\
q_{21} = -\displaystyle
\frac{x^1 x^2 y^2}{\Delta_1^2 \Delta_2^2}, \;\;\;
 p_{22} =  \displaystyle \frac{1}{\Delta_1 \Delta_2}, \;\;\;
q_{22} = \displaystyle \frac{x^2 y^2}{\Delta_1^2 \Delta_2^2},
      \end{array}
\right.
 \renewcommand{\arraystretch}{1}
\end{equation}
\begin{equation}\label{eq:180}
\renewcommand{\arraystretch}{1.3}
\left\{
\begin{array}{ll}
p = \displaystyle
\frac{1}{2 \Delta_1^2 \Delta_2^2} \Bigl[(x^1 - y^1) (1 + y^2) - x^1 (x^2  + y^1)\Delta_1
+ x^2 y^1 y^2)\Bigr], \\
 q = \displaystyle
\frac{1}{2 \Delta_1^2 \Delta_2^2} \Bigl[x^1\Bigl(1 + y^2 - (y^2)^2(x^2
+1)\Bigr) + y^1 \Bigl(1+ 2 y^2\Bigr) \Bigr],
        \end{array}
        \right.
\renewcommand{\arraystretch}{1}
\end{equation}
\begin{equation}\label{eq:181}
\renewcommand{\arraystretch}{1.3}
\left\{
\begin{array}{ll}
b^1_{111} = \displaystyle \frac{1}{2 \Delta_1^2\Delta_2^2}
\Bigl[(x^1 y^1 + x^2 y^2)(2 y^2 + x^1 (1 + y^2))\Bigr], \\
b^2_{111} = \displaystyle \frac{1}{2 \Delta_1^2\Delta_2^2}
\Bigl[y^1 y^2 (x^1 - 3 x^2 - 2 (x^1)^2)
+ 2 x^1 x^2 y^2 (1 - y^2) - (x^2)^2  y^1\Bigr], \\
b^1_{112} = \displaystyle \frac{1}{2 \Delta_1^2 \Delta_2^2}
\Bigl[x^2 y^2 (x^1 y^1 + x^2 y^2 - 3 x^1 y^2)
+ x^1 y^2 (y^2 - 2 x^1 y^1)\Bigr], \\
b^1_{121} = \displaystyle \frac{x^1 (x^2 - 1)}{2 \Delta_1
\Delta_2^2}, \;\;
b^1_{211} = \displaystyle \frac{x^2 (y^2)^2 - x^1 \Delta_1}{2
\Delta_1^2 \Delta_2^2},\;\;
b^2_{222} = \displaystyle \frac{1 - y^2}{(y^2)^2 \Delta_2^2},\\
b^1_{122} = b^1_{212} = -\displaystyle \frac{1}{2 \Delta_1 \Delta_2}, \;\;
   b^1_{22i} = b^2_{221} = 0,\\
b^2_{112} = \displaystyle \frac{1}{2 \Delta_1^2 \Delta_2^2}
\Bigl[x^1 y^2 (2 x^2 y^2 - 2 x^1 - x^2 - y^2 - 1
+  x^2 y^1 - 3  x^1 y^2) \\
\hspace*{24mm} + y^1 y^2 (x^2 - 4 x^1) - (x^1)^2 y^1 +
x^2 y^2 (x^2 y^2 - y^2 - 1)\Bigr], \\
b^2_{121} = \displaystyle \frac{1}{2 \Delta_1 \Delta_2^2 y^2}
\Bigl[x^1 (2 y^1 - 2 x^1 (1 + y^2) - y^2) +
 x^2 (x^2 +  y^2(3  -  x^1))\Bigr], \\
b^2_{211} = \displaystyle \frac{1}{2 \Delta_1^2 \Delta_2^2}
\Bigl[x^1 y^2 (y^2 (y^1 - x^1 - 2) - x^2 y^1 - 2) \\
\hspace*{24mm} +  x^2 (y^2 (x^2 + y^2 -  y^1) - (y^2)^2 + x^2 - y^1)\Bigr], \\
b^2_{122} = \displaystyle \frac{1}{2 \Delta_1 \Delta_2^2 y^2}
\Bigl[x^1 (2  - y^1 - x^1 y^2 - x^2 y^2)
- y^1 (1 + y^2)\Bigr], \\
b^2_{212} =
\displaystyle \frac{1}{2\Delta_1^2 \Delta_2^2}
\Bigl[y^2 (x^1 (x^2 (1 + y^2) - 1 + y^2) + 2 y^1)\Bigr], \\
        \end{array}
        \right.
\renewcommand{\arraystretch}{1}
\end{equation}

Since $b^2_{111} \neq 0$ (see (181)), it follows that
the web (175) is nontransversally geodesic. So, it belongs to the
class {\bf C}.

It follows from (180) that
$$
 p - q = \displaystyle
\frac{1}{2 \Delta_1^2 \Delta_2^2} \Bigl[ x^2 y^1 y^2
- x^1 y^2 (1 + y^2)(x^2 + y^1) + x^1 (y^2)^2 (1 + x^2)
- y^1 (3 y^2 + 2)\Bigr],
 $$
i.e., $p \neq q$, and  the inequalities
$(42)$ are satisfied.

It follows from (36)  that
$$
d p = p (\omega_1^1  + \omega_2^2) + \crazy{p}{1}_i^{ } \crazy{\omega}{1}^i +
 \crazy{p}{2}_i^{ } \crazy{\omega}{2}^i ,\;\;\;
 d q = q (\omega_1^1 + \omega_2^2) + \crazy{q}{1}_i^{ } \crazy{\omega}{1}^i +
 \crazy{q}{2}_i^{ } \crazy{\omega}{2}^i,
$$
where  $\omega_1^1$ and $\omega_2^2$ are defined by (177).
On the other hand, we can find $d p$ and $d q$ by
differentiating their expressions (180).
Comparing the results, we find that
\begin{equation}\label{eq:182}
 \renewcommand{\arraystretch}{1.3}
\left\{
\begin{array}{ll}
\crazy{p}{1}_1 =  \displaystyle\frac{1 - (x^2 + y^1) y^2 + x^1 y^1 -
\displaystyle\frac{(y^1)^2 y^2}{\Delta_1}}{2\Delta_1^2 \Delta_2^2}
+ \displaystyle\frac{p (3 y^1 + 5 x^1 y^2)}{\Delta_1 \Delta_2},
\\
\crazy{p}{1}_2 =  \displaystyle\frac{y^1 - x^1 (1 + y^2)}{2\Delta_1^2 \Delta_2^2}
- \displaystyle\frac{3 p}{y^2 \Delta_2},
\\
\crazy{q}{1}_1 =  \displaystyle\frac{1}{2\Delta_1^2
\Delta_2^3} \Bigl[x^1 (- x^1 + y^1 + x^1 (x^2 + y^1)(2y^2 + 1) - x^2 (y^1 + y^2 + 1))\\
\hspace*{22mm}  - x^2 (x^2 y^2 + y^2 + 1)\Bigr]
+ \displaystyle\frac{3 p x^1 (1 + 2 y^2)}{\Delta_1 \Delta_2},
\\
\crazy{q}{1}_2 =  \displaystyle\frac{1}{2\Delta_1^2
\Delta_2^3} \Bigl[x^1 (2 - (x^2 + y^1)(2 y^2 + 1) +
x^1 \Delta_1 + y^2 (1 - x^2) \\
\hspace*{22mm}- y^1 (1 + x^2) \Bigr]
- \displaystyle\frac{p (3 + 5 y^2)}{\Delta_1 \Delta_2},
\\
\crazy{p}{2}_1 =  \displaystyle\frac{y^2}{2 \Delta_1^3
\Delta_2^2} \Bigl[1 + y^2 (1 - y^2 (x^2 + 1) + x^1 y^1 y^2\bigr)\Bigr]
+ \displaystyle\frac{3 q (y^1 + 2 x^1 y^2)}{\Delta_1 \Delta_2},
\\
\crazy{p}{2}_2 = - \displaystyle\frac{x^1 y^2}{2\Delta_1^2 \Delta_2^2}
- \displaystyle\frac{3 q}{\Delta_1},
\\
\crazy{q}{2}_1 =  \displaystyle\frac{1}{2\Delta_1^2
\Delta_2^3} \Bigl[(x^1)^2 (2y^2 (x^2 + 1) - 1\bigr)
 + x^2 (1 + 2 y^2) - 2 x^1 y^1\Bigr] \\
\hspace*{7mm} + \displaystyle\frac{3 q x^1 (1 + 2 y^2)}{\Delta_1 \Delta_2},
\\
\crazy{q}{2}_2 =  \displaystyle\frac{1}{2\Delta_1^2
\Delta_2^3} \Bigl[x^1 \bigl(1 - 2 y^2 (x^2 + 1) - 2 y^2 - 1\bigr)
 + 2 y^1\Bigr] + \displaystyle\frac{2 q [1 +  y^2(3 - x^2)]}{\Delta_1
 \Delta_2}.
     \end{array}
\right.
 \renewcommand{\arraystretch}{1}
\end{equation}

If we substitute the values of $a_i$ from (178), $p,\, q$ from (180),
and  $\crazy{p}{\alpha}_1,\, \crazy{q}{\alpha}_1,\, \alpha = 1, 2,$
from (182) into the first equation (43), multiply the result by
the common denominator, and collect similar terms,
 we will see that there is  the term with $12(x^1)^4 \neq 0$.
 Thus the web (175) belongs to the class ${\bf G_4}$.

It is easy to prove by means of equations (178) and (179) that
conditions (8) do not hold. In fact, if we substitute the values
of $a_i$ from (178) and $p_{ij}$ from (179) into the first
 equation of (8), we can observe that there is a term $(x^2 y^2)^4$
which do not have similar terms. Thus the web
(175) belongs to the class ${\bf B}$.
}

We will present the results of this section in the following table
in which  we indicate to which classes the webs of our 18  examples
belong.

\vspace*{2mm}

\begin{center}
\begin{tabular}{|c|c|l|l|l|l|l|l|}
  \hline
  Example/Class & {\bf A} & {\bf B} & {\bf C} & {\bf D} & {\bf E} & {\bf F} & {\bf G} \\ \hline \hline
  1 &   &{\bf B}   & {\bf C} &  & ${\bf E_{11}}$ &  {\bf F} &  \\ \hline
  2 & ${\bf A_{21} \cap  A_{22}}$    &  & {\bf C} &  & ${\bf E_{22}}$  &  & ${\bf G_1}$ \\ \hline
  3 &    & {\bf B} & {\bf C} &  &  ${\bf E_{12}}$ &  & ${\bf G_1}$ \\ \hline
  4 & ${\bf A_{21} \cap A_{22}}$   &  &{\bf C}  &  &  ${\bf E_{22}}$ &  & ${\bf G_2}$ \\ \hline
  5 & ${\bf A_{31} \cap  A_{32}}$   &  &{\bf C} &  &  ${\bf E_{31}}$ & &  ${\bf G_3}$  \\ \hline
  6 & ${\bf A_{21} \cap  A_{22}}$  &  & {\bf C} &  &  ${\bf E_{23}}$ &  & ${\bf G_3}$  \\ \hline
  7 & ${\bf A_{31} \cap  A_{32}}$  &  & {\bf C} &  &  ${\bf E_{32}}$ &  & ${\bf G_3}$  \\ \hline
  8 &          &  & ? &? &   &  &  ${\bf G_3}$ \\ \hline
  9 &  ${\bf A_{21}}$ &  &?  &?  &  ${\bf E_{23}}$ &    &  ${\bf G_3}$ \\ \hline
  10 &  & {\bf B} & {\bf C}&  &  ${\bf E_{131}}$ &  &  ${\bf G_3}$  \\ \hline
  11 & &  {\bf B} & {\bf C} & &  ${\bf E_{111}}$ &  &  ${\bf G_3}$ \\ \hline
  12 &${\bf A_{131} \cap A_{132}}$ &  & {\bf C} &  &    &  &  ${\bf G_3}$ \\ \hline
  13 & ${\bf A_{31} \cap A_{32}}$ &  & {\bf C} &  &  ${\bf E_{32}}$ &  & ${\bf G_3}$  \\ \hline
   14 &  ${\bf A_{31} \cap A_{32}}$ &  &{\bf C} & & ${\bf E_{321}}$  &  & ${\bf G_3}$ \\ \hline
  15 & ${\bf A_{31} \cap A_{32}}$ &  & {\bf C} & &  ${\bf E_{32}}$ &  &  ${\bf G_3}$ \\ \hline
   16 & ${\bf A_{31}}$ &  &${\bf C}$ & &  ${\bf E_{3}}$ &  & ${\bf G_4}$ \\ \hline
   17 &  &  {\bf B}&{\bf C} & &${\bf E_{13}}$  &  & ${\bf G_4}$ \\ \hline
   18 &  & ${\bf B}$  & {\bf C} & &   &  &  ${\bf G_4}$ \\ \hline
  \hline
\end{tabular}
\end{center}

\vspace*{2mm}

Thus the examples considered in this section
 prove the existence of all webs indicated in this table.
Moreover, this proves the  existence of more general webs than
those indicated in the table. For example, the existence of
${\bf A_{131}}$ proves the existence of
${\bf A_{13}}$ and  ${\bf A_{1}}$.

\setcounter{theorem}{14}

\rem{  We can see from the
 above table that the 3-webs
of Examples {\bf 7} and {\bf 13} belong to the same classes.
In order to show that they are not equivalent,
we compare the vanishing components of their
curvature tensor. In Example {\bf 7}, they are
$$
b^2_{jkl} = b^1_{1kl}= b^1_{211} = 0.
$$
In Example {\bf 13}, they are
$$
b^2_{jkl} = b^1_{11i} =  b^1_{11i} = 0.
$$
Thus the 3-webs in these two examples are not equivalent.

We also can see from the above table that the 3-webs
of Examples {\bf 13} and {\bf 15} belong to the same classes.
It follows from (141) and (155) that they
have the same vanishing  components of their
curvature tensors. So, some additional investigation is needed
in order to determine whether  the 3-webs
of Examples {\bf 13} and {\bf 15} are equivalent.
}

\rem{In the table above for the webs of Examples 8 and 9,
 we put the question marks in the columns {\bf C} and {\bf D}
 since for certain values of the coefficients
 $c^i_{jk}$ the polynomial webs of these examples belong
 to the class {\bf C}, and for other values of these
 coefficients, they belong to the class {\bf D}.

}

One can also see from the table above
that in our examples there is no 3-webs
that are  transversally geodesic (Class ${\bf D}$).
This gives a rise to the following problem:
{\em Construct an example of a nonisoclinic transversally
geodesic nonhexagonal $($or hexagonal$)$ $3$-web.}

\noindent {\em Author's address}:

\noindent Vladislav V. Goldberg\\ Department of Mathematics\\ New
Jersey Institute of Technology\\ Newark, N.J. 07102, U.S.A.

\vspace*{2mm}

\noindent {\em Author's e-mail address}: vlgold@m.njit.edu

\end{document}